\documentclass{amsart}
\usepackage{amsthm}
\usepackage{amssymb}
\usepackage{upref}
\usepackage{amscd}

\numberwithin{equation}{section}
\newtheorem{main}{Theorem}

\newtheorem{theorem}[equation]{Theorem}
\newtheorem{lemma}[equation]{Lemma}
\newtheorem{proposition}[equation]{Proposition}
\newtheorem{corollary}[equation]{Corollary}
\newtheorem{conjecture}[equation]{Conjecture}
\theoremstyle{definition}
\newtheorem{definition}[equation]{Definition}
\newtheorem*{remark}{Remark}
\newtheorem{example}[equation]{Example}

\DeclareMathOperator{\Hom}{Hom}
\DeclareMathOperator{\Ann}{Ann}
\DeclareMathOperator{\Prim}{Prim}

\DeclareMathOperator{\Ext}{Ext}
\DeclareMathOperator{\Spec}{Spec}

\DeclareMathOperator{\colim}{colim}

\DeclareMathOperator{\im}{im}

\DeclareMathOperator{\coker}{coker}
\DeclareMathOperator{\height}{ht}

\renewcommand{\smash}{\wedge}

\newcommand{\Q}{\mathbb{Q}}

\renewcommand{\ell}{Ell}

\newcommand{\cat}[1]{\mathcal{#1}}
\newcommand{\Ring}{\mathbf{Rings}}

\newcommand{\Aff}{\mathbf{Aff}}

\newcommand{\Z}{\mathbb{Z}}
\newcommand{\Fp}{\mathbb{F}_{p}}
\newcommand{\Zp}{\Z _{(p)}}
\newcommand{\comod}{\text{-comod}}

\newcommand{\mathcolon}{\colon\,}
\newcommand{\suchthat}{\ \vert \ }

\newcommand{\usc}{\textup{;}}

\begin{document}
 
\title{Comodules and Landweber exact homology theories} 

\date{\today}

\author{Mark Hovey}
\address{Department of Mathematics \\ Wesleyan University
\\ Middletown, CT 06459}
\email{hovey@member.ams.org}

\author{Neil Strickland}
\address{Department of Pure Mathematics \\ University of Sheffield
\\ Sheffield S3 7RH \\
England}
\email{N.P.Strickland@sheffield.ac.uk}


\begin{abstract}
We show that, if $E$ is a commutative $MU$-algebra spectrum such that
$E_{*}$ is Landweber exact over $MU_{*}$, then the category of
$E_{*}E$-comodules is equivalent to a localization of the category of
$MU_{*}MU$-comodules.  This localization depends only on the heights of
$E$ at the integer primes $p$.  It follows, for example, that the
category of $E(n)_{*}E(n)$-comodules is equivalent to the category of
$(v_{n}^{-1}BP)_{*}(v_{n}^{-1}BP)$-comodules.  These equivalences give
simple proofs and generalizations of the Miller-Ravenel and Morava
change of rings theorems.  We also deduce structural results about the
category of $E_{*}E$-comodules.  We prove that every $E_{*}E$-comodule
has a primitive, we give a classification of invariant prime ideals in
$E_{*}$, and we give a version of the Landweber filtration theorem.
\end{abstract}

\maketitle

\section*{Introduction}

Suppose $E_{*}(-)$ and $R_{*}(-)$ are reduced homology theories with
commutative products defined on finite CW complexes.  Then $E_{*}(-)$ is
said to be \textbf{Landweber exact} over $R_{*}(-)$ if there is a natural
isomorphism
\[
E_{*}(X) \cong E_{*}\otimes _{R_{*}} R_{*}(X)
\]
for all finite CW complexes $X$.  It then follows that this natural
isomorphism extends to all CW complexes $X$, and indeed to all spectra
$X$.  Because of this, we usually just say that the spectrum $E$ is
Landweber exact over the spectrum $R$.  Examples of this phenomenon
abound in stable homotopy theory, and were first studied in some
generality by Landweber~\cite{land-exact}.  

\begin{example}\label{ex-exactness}
\begin{enumerate}
\item [(a)] Conner and Floyd~\cite{conner-floyd} showed that complex
$K$-theory $K$ is Landweber exact over complex cobordism $MU$, and also
that real $K$-theory $KO$ is Landweber exact over symplectic cobordism
$MSp$.  Hopkins and the first author~\cite{hovey-Spin} showed that $KO$
is also Landweber exact over Spin cobordism $MSpin$.
\item [(b)] The various elliptic cohomology
theories~\cite{landweber-ravenel-stong} are all Landweber exact over
$MU$.  
\item [(c)] The Brown-Peterson spectrum $BP$ at a prime $p$ is Landweber
exact over $MU$.  Furthermore, the $p$-localization $MU_{(p)}$ of $MU$
is also Landweber exact over $BP$~\cite{land-exact}.
\item [(d)] The Johnson-Wilson spectrum $E(n)$~\cite{johnson-wilson-en}
as well as the Morava $E$-theory spectrum $E_{n}$ used in the work of
Hopkins and Miller~\cite{hopkins-miller} are Landweber exact over $BP$. 
\item [(e)] The Morava $K$-theory spectrum $K(n)$ is Landweber exact
over the spectrum $P(n)=BPI_{n}$~\cite{yosimura}.  
\end{enumerate}
\end{example}

In all of the examples of spectra $E$ above, the module $E_*E$ is flat
over $E_*$, so $(E_*,E_*E)$ is a Hopf algebroid, or equivalently, a
groupoid object in the opposite of the category of graded-commutative
rings.  For compatibility with the usual conventions in topology, we
set up this correspondence so that the maps $\eta_L,\eta_R\mathcolon
E_*\xrightarrow{}E_*E$ represent the maps sending a morphism to its
target and source respectively.  We refer to~\cite[Appendix
1]{ravenel} for basic facts about Hopf algebroids.  The reduced
homology $E_{*}X$ is a comodule over the Hopf algebroid
$(E_{*},E_{*}E)$~\cite[Proposition~2.2.8]{ravenel}.  One of the main
reasons this is important is because the $E_{2}$-term of the Adams
spectral sequence of $X$ based on $E$ is
\[
\Ext_{E_{*}E}^{**}(E_{*}, E_{*}X),
\]
and this $\Ext $ is taken in the category of $E_{*}E$-comodules.   
To help compute these $E_2$-terms, various authors have constructed
isomorphisms of the form
$\Ext_\Gamma^{**}(M,N)\simeq\Ext_\Sigma^{**}(M',N')$ under various
hypotheses on the algebroids $\Gamma$ and $\Sigma$, and the comodules
$M$, $N$, $M'$ and $N'$.  This includes the change of rings theorems
of Miller-Ravenel~\cite{miller-ravenel}, Morava~\cite{morava-comodules},
and the first author and Sadofsky~\cite{hovey-sadofsky-picard}.  The
main result of this paper is that these isomorphisms come from
equivalences of comodule categories, and that such equivalences are
much more common and systematic than was previously suspected.

The definition of Landweber exactness given above for homology
theories has an analogue for Hopf algebroids.  Given a Hopf algebroid
$(A, \Gamma)$ and an $A$-algebra $B$, we define $B$ to be
\textbf{Landweber exact} over $(A, \Gamma)$, or, by abuse of notation,
over $A$, if the functor $B\otimes_{A}(-)$ is exact on the category of
$\Gamma$-comodules.  If $E_{*}(-)$ and $R_{*}(-)$ are homology
theories as above, $E$ is an $R$-module spectrum, and $E_{*}$ is
Landweber exact over $R_{*}$, then a well-known argument shows that
$E$ is Landweber exact over $R$.  In the examples listed above,
$K_{*}$, $\ell_{*}$, and $BP_{*}$ are all Landweber exact over
$MU_{*}$, and $(MU_{*})_{(p)}$, $E(n)_{*}$, and $E_{n*}$ are Landweber
exact over $BP_{*}$, and $K(n)_{*}$ is Landweber exact over
$P(n)_{*}$.  On the other hand, it is not known whether $KO_{*}$ is
Landweber exact over $MSp_{*}$ or $MSpin_{*}$.

In the above situation, it is important to understand the relationship
between $E_*E$-comodules and $R_*R$-comodules.  Given a Hopf algebroid
$(A,\Gamma)$ and a ring map $f\mathcolon A\xrightarrow{}B$, we put
$\Gamma_B=B\otimes_A\Gamma\otimes_AB$.  The pair $(B,\Gamma_B)$ has a
natural structure as a Hopf algebroid; we recall some details in
Section~\ref{sec-landweber}.  A central object of this paper is to
make a detailed study of the relationship between the category of
$\Gamma_{B}$-comodules and the category of $\Gamma$-comodules when $B$
is Landweber exact over $A$.  We prove the following theorem as
Theorem~\ref{thm-Giraud}.

\begin{main}\label{main-A}
Suppose $(A, \Gamma)$ is a flat Hopf algebroid and $B$ is a Landweber
exact $A$-algebra.  Then the category of $\Gamma_{B}$-comodules is
equivalent to the localization of the category of $\Gamma$-comodules
with respect to the hereditary torsion theory 
\[
\cat{T}=\{M\suchthat B\otimes_{A}M=0 \}.
\]
\end{main}

To apply this to cases of interest in algebraic topology, we give a
partial classification of graded hereditary torsion theories of
$BP_{*}BP$-comodules in Theorem~\ref{thm-torsion}. 

\begin{main}\label{main-B}
 Let $\cat{T}$ be a graded hereditary torsion theory of
 $BP_{*}BP$-comodules, and suppose that $\cat{T}$ contains a nonzero
 finitely presented comodule.  Then either $\cat{T}$ is the whole
 category of comodules, or there is an $n$ such that $\cat{T}$ is the
 collection of $v_{n}$-torsion comodules.
\end{main}

This theorem is analogous to the classification of thick subcategories
of finite $p$-local spectra~\cite{hopkins-smith}.

This then leads to our main result, which is Theorem~\ref{thm-main-BP}.  

\begin{main}\label{main-C}
Define the \textbf{height} of a Landweber exact $BP_{*}$-algebra $E_{*}$
to be the largest $n$ such that $E_{*}/I_{n}$ is nonzero.  If $E_{*}$
and $E'_{*}$ are Landweber exact $BP_{*}$-algebras of the same height,
then the category of graded $E_{*}E$-comodules is equivalent to the
category of graded $E'_{*}E'$-comodules.  In particular, the categories
of $E(n)_{*}E(n)$-comodules, $E_{n*}E_{n}$-comodules, and
$(v_{n}^{-1}BP)_{*}(v_{n}^{-1}BP)$-comodules are all equivalent.  
\end{main}

As mentioned previously, this gives a simple explanation for the
change of rings theorems of Miller-Ravenel, Morava, and
Hovey-Sadofsky, all of which say that two $\Ext$ groups computed over
different Hopf algebroids are isomorphic.  Namely, the $\Ext$ groups
are isomorphic because the categories they are computed in are
equivalent.

When $E_{*}$ is Landweber exact over $BP_{*}$, the category of
$E_{*}E$-comodules is a localization of the category of
$BP_{*}BP$-comodules, by Theorem~\ref{main-A}.  This allows us to extend
the standard structure theorems for $BP_{*}BP$-comodules of
Landweber~\cite{land-exact} to $E_{*}E$-comodules.  The following
theorem is a summary of the results of Section~\ref{sec-structure}.  

\begin{main}\label{main-D}
Suppose $E_{*}$ is a Landweber exact $BP_{*}$-algebra of height $n\leq \infty$.
\begin{enumerate}
\item [(a)] Every nonzero $E_{*}E$-comodule has a nonzero primitive. 
\item [(b)] If $I$ is an invariant radical ideal in $E_{*}$, then
$I=I_{k}$ for some $k\leq n$.  
\item [(c)] Every $E_{*}E$-comodule $M$ that is finitely presented over
$E_{*}$ admits a finite filtration by subcomodules 
\[
0=M_{0} \subseteq M_{1} \subseteq \dotsb \subseteq M_{s} = M
\]
for some $s$, with $M_{r}/M_{r-1}\cong s^{t}E_{*}/I_{j}$ for some $j\leq
n$ and some $t$, both depending on $r$.  
\end{enumerate}
\end{main}

\begin{remark}
 Baker~\cite{baker-comodules} has constructed a counterexample to a
 statement closely related to~(a), in the case where $E$ is the Morava
 $E$-theory spectrum $E_n$.  This is not in fact a contradiction,
 because of the difference between $E_*E=\pi_*(E\smash E)$ (which is
 used in our work) and $\pi_*L_{K(n)}(E\smash E)$ (which is more
 closely related to the Morava stabiliser group, and is used in
 Baker's work).  The topological comodule categories considered by
 Baker are well-known to be important, but they do not fit into our
 present framework; we hope to return to this in future.
\end{remark}

The theorems we have just discussed all have global versions, where we
replace $BP_{*}$ by $MU_{*}$, and more local versions, where we replace
$BP_{*}$ by $BP_{*}/J$ for a nice regular sequence $J$.  We discuss
these versions briefly at the end of the paper.  

As mentioned above, the category of $E(n)_{*}E(n)$-comodules is a
localization of the category of $BP_{*}BP$-comodules.  The resulting
localization functor on $BP_{*}BP$-comodules is denoted $L_{n}$, and
is analogous to the chromatic localization functor $L_{n}$ much used
in stable homotopy theory~\cite{rav-nilp}.  The algebraic $L_{n}$ is
very interesting in its own right; it is left exact, and has
interesting right derived functors $L_{n}^{i}$, which are closely
related to local cohomology.  The functor $L_{n}$ and its derived
functors are studied in~\cite{hovey-strickland-derived}.

We also point out that, to give a good algebraic model for stable
homotopy theory, one wants a triangulated category rather than an abelian
category.  So there should be analogues of the theorems in this paper
for some kind of derived categories of $BP_{*}BP$-comodules and
$E_{*}E$-comodules.  There are problems with the ordinary derived
category; the first author has constructed a well-behaved replacement
for it in~\cite{hovey-comodule-homotopy}.  The authors have proved
analogues of some of the theorems of this paper for these derived
categories in~\cite{hovey-barcelona}.  

The authors would like to thank the Universitat de Barcelona, the
Universitat Aut\`{o}noma de Barcelona, the Centre de Recerca Matematica,
and the Isaac Newton Insitute for Mathematical Sciences for their
support during this project.  They would also like to thank John
Greenlees and Haynes Miller for several helpful discussions about this
paper.  

\section{Abelian localization}\label{sec-local}

In this section, we summarize Gabriel's theory of localization of
abelian categories from an algebraic topologist's point of view for the
convenience of the reader.  The original source for this material
is~\cite{gabriel-abelian-categories}; a standard source for localization
in module categories is~\cite{stenstrom}.  The
book~\cite{van-oyst-localization} gives a quick summary of the theory in
an arbitrary Grothendieck category.  

The following definition is standard in homotopy theory. 

\begin{definition}\label{defn-local}
Suppose $\cat{E}$ is a class of maps in a category $\cat{C}$.  An object
$X$ of $\cat{C}$ is said to be \textbf{$\cat{E}$-local} if
$\cat{C}(f,X)$ is an isomorphism of sets for all $f\in \cat{E}$.  We
denote the full subcategory of $\cat{E}$-local objects by
$L_{\cat{E}}\cat{C}$.  An \textbf{$\cat{E}$-localization} of an object
$M\in \cat{C}$ is a map $M\xrightarrow{}LM$ in $\cat{E}$ where $LM\in
L_{\cat{E}}\cat{C}$.  If every $M\in \cat{C}$ has an
$\cat{E}$-localization, we say that \textbf{$\cat{E}$-localizations
exist}.
\end{definition}

It is also possible to define localizations without reference to the
class $\cat{E}$.  

\begin{definition}\label{defn-localization-functor}
A \textbf{localization functor} on a category $\cat{C}$ is a
functor $L\mathcolon  \cat{C}\xrightarrow{}\cat{C}$ and a natural
transformation $\iota_{M} \mathcolon M\xrightarrow{}LM$ such that
$L\iota_{M}=\iota_{LM}$ and this map is an isomorphism.  
\end{definition}

The following proposition is reasonably well-known; a version of it can
be found in~\cite[Section~3.1]{hovey-axiomatic} and in other places.

\begin{proposition}\label{prop-local}
Suppose $L$ is a localization functor on a category $\cat{C}$.  Let
$\cat{E}$ denote the class of all maps $f$ such that $Lf$ is an
isomorphism.  Then $\iota_{M}$ is an $\cat{E}$-localization of $M$ for
all $M\in \cat{C}$.  Conversely, if $\cat{E}$ is a class of maps on
$\cat{C}$ such that an $\cat{E}$-localization $\iota_{M}\mathcolon
M\xrightarrow{}LM$ exists for all $M\in \cat{C}$, then $L$ is a
localization functor.  Furthermore, in either case $L$, thought of as a
functor $L\mathcolon \cat{C}\xrightarrow{}L_{\cat{E}}\cat{C}$, is left
adjoint to the inclusion functor.
\end{proposition}

We refer to the localization functor of Proposition~\ref{prop-local} as
localization \textbf{with respect to $\cat{E}$}.  

A common way for localizations to arise is displayed in the following
proposition.  

\begin{proposition}\label{prop-local-adjoint}
Suppose $F\mathcolon \cat{C}\xrightarrow{}\cat{D}$ is a functor with
right adjoint $G$, and the counit of the adjunction
$\epsilon_{M}\mathcolon FGN\xrightarrow{}N$ is an isomorphism for all
$M\in \cat{D}$.  Then $GF$ is the localization functor on $\cat{C}$ with
respect to $\cat{E}=\{f| Ff \text{ is an isomorphism} \}$.  Furthermore,
$G$ defines an equivalence of categories $G\mathcolon
\cat{D}\xrightarrow{}L_{\cat{E}}\cat{C}$.
\end{proposition}

\begin{proof}
Let $L=GF$. The natural transformation $\iota_{M}\mathcolon
M\xrightarrow{}LM$ is the unit $\eta_{M}$ of the adjunction.  The two
triangular relations of the adjunction say, respectively, that
\[
\epsilon_{FM}\circ (F\eta_{M})=1_{FM} \text{ and } (G\epsilon_{M})\circ
\eta_{GM}=1_{GM}. 
\]
In particular, $GF\eta_{M}=(G\epsilon_{FM})^{-1}=\eta_{GFM}$.  This
means that $L\iota_{M}=\iota_{LM}$ and that this map is an isomorphism,
as required.  By Proposition~\ref{prop-local}, $L$ is localization with
respect to $\cat{E}=\{f|Lf \text{ is an isomorphism} \}$.  Since $FG$ is
naturally isomorphic to the identity, one can easily check that $GFf$ is
an isomorphism if and only if $Ff$ is an isomorphism.  

Since $FG$ is naturally isomorphic to the identity, $G$ defines an
equivalence of categories from $\cat{D}$ to its image.  Adjointness
shows that $GN$ is $\cat{E}$-local for all $N\in \cat{D}$.  Conversely,
the image of $G$ contains $LM$ for all $M\in \cat{C}$, so is a skeleton
of $L_{\cat{E}}\cat{C}$.  The result follows.  
\end{proof}

In point of fact, every localization functor arises from an adjunction
as in Proposition~\ref{prop-local-adjoint}; if $L$ is a localization
functor on $\cat{C}$, we can think of it as a functor $L\mathcolon
\cat{C}\xrightarrow{}L_{\cat{E}}\cat{C}$, where it is left adjoint to
the inclusion and satisfies the hypotheses of
Proposition~\ref{prop-local-adjoint}.  

Now suppose that our category $\cat{C}$ is abelian.  It is natural,
then, to consider localization functors arising from adjunctions
$F\mathcolon \cat{C}\rightleftarrows\cat{D}\mathcolon G$ as in
Proposition~\ref{prop-local-adjoint} where $\cat{D}$ is also abelian
and $F$ is exact.

\begin{definition}\label{defn-torsion}
 Suppose $\cat{T}$ is a full subcategory of an abelian category
 $\cat{C}$.  Then $\cat{T}$ is called a \textbf{hereditary torsion
   theory} if $\cat{T}$ is closed under subobjects, quotient objects,
 extensions, and arbitrary coproducts.  When $\cat{T}$ is a hereditary
 torsion theory, we define the class $\cat{E}_{\cat{T}}$ of
 \textbf{$\cat{T}$-equivalences} to consist of those maps whose
 cokernel and kernel are in $\cat{T}$.  We define an object to be
 \textbf{$\cat{T}$-local} if and only if it is
 $\cat{E}_{\cat{T}}$-local.  We let $L_{\cat{T}}\cat{C}$ denote the
 full subcategory of $\cat{T}$-local objects.
\end{definition}

Note that a hereditary torsion theory is just a Serre class that is
closed under coproducts.  Also, one can form the smallest hereditary
torsion theory containing a specified class of objects by taking the
intersection of all hereditary torsion theories containing that class.

\begin{proposition}\label{prop-local-torsion}
Suppose $\cat{C}$ and $\cat{D}$ are abelian categories, $F\mathcolon
\cat{C}\xrightarrow{}\cat{D}$ is an exact functor with right adjoint
$G$, and the counit of the adjunction $\epsilon_{M}\mathcolon
FGN\xrightarrow{}N$ is an isomorphism for all $M\in \cat{D}$.  Then $GF$
is the localization functor on $\cat{C}$ with respect to the hereditary
torsion theory $\cat{T}=\ker F =\{M|FM=0 \}$.  Furthermore, $G$ defines
an equivalence of categories $G\mathcolon
\cat{D}\xrightarrow{}L_{\cat{T}}\cat{C}$.
\end{proposition}

\begin{proof}
Proposition~\ref{prop-local-adjoint} implies that $GF$ is localization
with respect to 
\[
\cat{E}=\{f|Ff \text{ is an isomorphism} \}.
\]
But, since $F$ is exact, $Ff$ is an isomorphism if and only if $F(\ker
f)=F(\coker f)=0$, which is true if and only if $f$ is a
$\cat{T}$-equivalence.
\end{proof}

The main result of Gabriel on abelian localizations is the following
theorem.  Recall that a Grothendieck category is a cocomplete abelian
category with a generator in which filtered colimits are exact.

\begin{theorem}\label{thm-Gabriel}
Suppose $\cat{T}$ is a hereditary torsion theory in a Grothendieck
abelian category $\cat{C}$.  Then $\cat{T}$-localizations exist.  
\end{theorem}

We outline the proof of Gabriel's theorem~\ref{thm-Gabriel}, as we will
need some of the ideas from this proof later.  We first recall the
characterization of $\cat{T}$-local objects.

\begin{lemma}\label{lem-E-local}
Suppose $\cat{T}$ is a hereditary torsion theory in an abelian category
$\cat{C}$.  An object $X$ of $\cat{C}$ is $\cat{T}$-local if and only if
$\cat{C}(T,X)=\Ext ^{1}_{\cat{C}}(T,X)=0$ for all $T\in \cat{T}$.
\end{lemma}

Recall that one can define $\Ext $ in an arbitrary abelian category
without recourse to either projectives or
injectives~\cite{maclane-homology}. In particular,
$\Ext^{1}_{\cat{C}}(M,N)$ is just the collection of all equivalence
classes of short exact sequences
\[
0 \xrightarrow{} N \xrightarrow{} E \xrightarrow{} M \xrightarrow{} 0.
\]
The usual exact sequences for $\Ext $ work in this generality. 

\begin{proof}
Suppose first that $X$ is $\cat{T}$-local, and $T\in \cat{T}$.  Since
$0\xrightarrow{}T$ is a $\cat{T}$-equivalence, we conclude that
$\cat{C}(T,X)=0$. 
Given a short exact sequence 
\[
0 \xrightarrow{} X \xrightarrow{f} Y \xrightarrow{} T \xrightarrow{}0,
\]
the map $f$ is a $\cat{T}$-equivalence, so
$f^*\mathcolon\cat{C}(Y,X)\xrightarrow{}\cat{C}(X,X)$ is an
isomorphism.  Thus the identity map of $X$ comes from a map
$g\mathcolon Y\xrightarrow{}X$, and $g$ defines a splitting of the
given sequence.  Hence $\Ext ^{1}_{\cat{C}}(T,X)=0$.

Conversely, suppose $\cat{C}(T,X)=\Ext ^{1}_{\cat{C}}(T,X)=0$ for all
$T\in \cat{T}$, and $f\mathcolon A\xrightarrow{}B$ is a
$\cat{T}$-equivalence.  Consider the two short exact sequences 
\[
0 \xrightarrow{} \ker f \xrightarrow{} A \xrightarrow{} \im f
\xrightarrow{} 0
\]
and
\[
0 \xrightarrow{} \im f \xrightarrow{} B \xrightarrow{} \coker f
\xrightarrow{} 0.
\]
By applying the functor $\cat{C}(-,X)$ to these short exact sequences,
using the fact that 
\[
\cat{C}(\ker f, X) =\cat{C}(\coker f,X) = \Ext^{1}_{\cat{C}}(\coker
f,X)=0,
\]
we see that $\cat{C}(f,X)$ is an isomorphism.  
\end{proof}

Now, in order to construct the localization $L_{\cat{T}}(X)$ of $X$ with
respect to a hereditary torsion theory $\cat{T}$ in a Grothendieck
category $\cat{C}$, we first form the union $TX$ of all the subobjects
of $X$ that are in $\cat{T}$ (these form a set because we are in a
Grothendieck category).  This gives us a $\cat{T}$-equivalence $X
\xrightarrow{}X/TX$.  Then we taken an injective envelope $I$ of $X/TX$
(injective envelopes exist in a Grothendieck category), 
producing an exact sequence 
\[
0 \xrightarrow{} X/TX \xrightarrow{} I \xrightarrow{} Q \xrightarrow{} 0.
\]
Finally, we let $L_{\cat{T}}(X)$ be the pullback $I\times _{Q} TQ$.  The
induced embedding $X/TX\xrightarrow{}L_{\cat{T}}(X)$ is a
$\cat{T}$-equivalence, and one can check that $L_{\cat{T}}(X)$ is
$\cat{T}$-local.  

\begin{remark}
 In our case we will be working with graded abelian categories
 $\cat{C}$.  This means that we have a given self-equivalence
 $s\mathcolon \cat{C}\xrightarrow{}\cat{C}$, which in fact is an
 isomorphism of categories in our examples.  In this case, we define a
 full subcategory $\cat{D}$ to be \textbf{graded} when $X\in \cat{D}$
 if and only if $sX\in \cat{D}$.  Similarly, a class of maps $\cat{E}$
 in $\cat{C}$ is said to be \textbf{graded} when $f\in \cat{E}$ if and
 only if $sf\in \cat{E}$.  The results of this section all have
 corresponding graded versions.
\end{remark}

\section{Landweber exact algebras}\label{sec-landweber}

In this section, we apply localization techniques to comodules over
Hopf algebroids.  Recall that a Hopf algebroid is a pair of (possibly
graded) commutative rings $(A,\Gamma)$ so that $\Ring (A, R)$ and
$\Ring (\Gamma ,R)$ are the objects and morphisms of a groupoid that
is natural in the (graded) commutative ring $R$.  The precise
structure maps and relations they satisfy can be found
in~\cite[Appendix~1]{ravenel}.  The reason we are interested in them
is that $(E_{*}, E_{*}E)$ is a Hopf algebroid for many homology
theories $E$, as explained in~\cite[Proposition~2.2.8]{ravenel}.

We will always assume our Hopf algebroids are \textbf{flat}; this
means that the left unit $\eta_{L}\mathcolon A\xrightarrow{}\Gamma$
corepresenting the target of a morphism is a flat ring extension.
This is the same as assuming that the right unit $\eta_{R}\mathcolon
A\xrightarrow{}\Gamma$ corepresenting the source of a morphism is
flat.

We note that in working with Hopf algebroids it is important to remember
that $M\otimes _{A}N$ always denotes the tensor product of
$A$-bimodules, using the right $A$-module structure on $M$ and the left
$A$-module structure on $N$.  This mostly matters for $\Gamma $, which
is a right $A$-module via the right unit $\eta _{R}$ and a left
$A$-module via the left unit $\eta _{L}$.  

Recall that a \textbf{comodule} over a Hopf algebroid $(A, \Gamma)$ is a
left $A$-module $M$ equipped with a coassociative and counital coaction
map $\psi \mathcolon M\xrightarrow{}\Gamma \otimes_{A}M$.  The category
of $\Gamma$-comodules is abelian when $(A, \Gamma)$ is
flat~\cite[Theorem~A1.1.3]{ravenel}.  

We now recall the definition of Landweber exactness, mentioned in the
introduction.
\begin{definition}\label{defn-Landweber}
 Suppose $(A,\Gamma)$ is a flat Hopf algebroid, and $f\mathcolon
 A\xrightarrow{}B$ is a ring homomorphism.  We say that $B$ is
 \textbf{Landweber exact} over $(A,\Gamma)$, or just over $A$, if the
 functor $M\mapsto B\otimes_{A}M$ from $\Gamma $-comodules to
 $B$-modules is exact.
\end{definition}

%



We next recall the construction of the Hopf algebroid $\Gamma_B$, and
use it to reformulate the notion of Landweber exactness.  The
definition is motivated by the following construction on groupoids.
Consider a groupoid with object set $X$ and morphism set $G$.  Given a
set $Y$ and a map $f\mathcolon Y\xrightarrow{}X$, we define a new
groupoid $(Y,G_f)$ as follows: the object set is $Y$, and the
morphisms in $G_f$ from $y_1$ to $y_0$ are the morphisms in $G$ from
$f(y_1)$ to $f(y_0)$, so as a set we have 
\[ G_f = Y \times_{X,f} G \times_{X,f} Y. \]  
The map $f$ induces a full and faithful functor
$F\mathcolon (Y,G_f)\xrightarrow{}(X,G)$.  To understand when this is an
equivalence, consider the set
\[ U =
   \{(y,g) \suchthat y\in Y\;,\; g\in G \;,\; f(y)=\text{target}(g) \}
   = Y\times_XG.
\]
There is a map $\pi\mathcolon U\xrightarrow{}X$ given by
$(y,g)\mapsto\text{source}(g)$.  Our functor $F$ is essentially
surjective, and thus an equivalence, iff $\pi$ is surjective.

Now suppose we have a Hopf algebroid $(A,\Gamma)$ and an $A$-algebra
$B$.  For any ring $T$, we have a groupoid
\[
(\Ring(A,T),\Ring(\Gamma,T))
\]
and a map $\Ring(B,T)\xrightarrow{}\Ring(A,T)$.  We can apply the
construction above to obtain a new groupoid
$(\Ring(B,T),\Ring(\Gamma_B,T))$, where
$\Gamma_B=B\otimes_A\Gamma\otimes_AB$ as before.  The groupoid structure
is natural in $T$, so Yoneda's lemma gives $(B,\Gamma_B)$ the structure
of a Hopf algebroid.  For further details,
see~\cite[p.~1315]{hovey-hopf}.  There is a morphism
$\Phi=(f,\tilde{f})\mathcolon(A,\Gamma)\xrightarrow{}(B,\Gamma_B)$ of
Hopf algebroids, where $\tilde{f}(u)=1\otimes u\otimes 1$; this
corresponds to the functor $F$.  The morphism $\Phi$ induces a functor
$\Phi_*$ from $\Gamma$-comodules to $\Gamma_B$-comodules, given by
$M\mapsto B\otimes_AM$; by definition, this is exact iff $B$ is
Landweber exact over $A$.

The map $\pi\mathcolon U\xrightarrow{}X$ corresponds to the ring map
$f\otimes\eta_R\mathcolon A\xrightarrow{}B\otimes_A\Gamma$ given by
$a\mapsto 1\otimes\eta_R(a)$.

\begin{lemma}\label{lem-land-exact}
 Suppose $(A,\Gamma )$ is a flat Hopf algebroid, and
 $f\mathcolon A\xrightarrow{}B$ is a ring homomorphism.  Then $B$ is
 Landweber exact if and only if the map $f\otimes\eta_R$ makes
 $B\otimes _{A}\Gamma $ into a flat $A$-algebra.
\end{lemma}

\begin{proof}
Suppose that $B$ is Landweber exact, and $M\xrightarrow{}N$ is a
monomorphism of $A$-modules.  Then $\Gamma \otimes_{A} f$ is a
monomorphism as well, since $\Gamma $ is flat as a right $A$-module.
But $\Gamma \otimes _{A}M$ and $\Gamma \otimes _{A}N$ are both $\Gamma
$-comodules, with the coaction coming from the diagonal on $\Gamma $.
(This is called the \textbf{extended comodule structure} on $\Gamma
\otimes _{A}M$).  This makes $\Gamma \otimes _{A}f$ a comodule map.
Since $B$ is Landweber exact, we conclude that $B\otimes _{A}\Gamma
\otimes _{A}f$ is a monomorphism.

Conversely, suppose that $B\otimes _{A}\Gamma $ is flat over $A$, and
$u\mathcolon M\xrightarrow{}N$ is a monomorphism of comodules.  The
coaction map $\psi _{M}$ is a split monomorphism of $A$-modules; the
splitting is given by $\epsilon \otimes 1$, where $\epsilon $ is the
counit of $(A,\Gamma )$.  Hence $u$ is a retract of $\Gamma \otimes_{A}
u$ as a map of $A$-modules.  It follows that $B\otimes_{A} u$ is a
retract of $B\otimes_{A} \Gamma \otimes_{A} u$ as a map of $B$-modules.
Since $B\otimes _{A}\Gamma$ is flat over $A$, we conclude that
$B\otimes _{A}u$ is a monomorphism, as required.
\end{proof}
\begin{corollary}
 If $B$ is Landweber exact over $A$, then $\Gamma_B$ is flat over
 $B$ (so the category of $\Gamma_B$-comodules is abelian). 
\end{corollary}
\begin{proof}
 We have seen that $B\otimes_A\Gamma$ is flat as a right $A$-module;
 now take tensor products with $B$ on the right.
\end{proof}



For any morphism $\Phi$ of flat Hopf algebroids, the functor
$\Phi_{*}$ obviously preserves colimits, so it should have a right
adjoint $\Phi^{*}$; we next check that this works.

\begin{lemma}\label{lem-adjoint}
Suppose $\Phi \mathcolon (A,\Gamma )\xrightarrow{}(B,\Sigma )$ is a map
of flat Hopf algebroids.  Then the functor $\Phi_{*}\mathcolon \Gamma
\comod \xrightarrow{}\Sigma \comod $ defined by
$\Phi_{*}M=B\otimes_{A}M$ has a right adjoint $\Phi^{*}$.
\end{lemma}

This lemma is proved in~\cite[Section~1]{hovey-comodule-homotopy}, but
it is central to our work, so we recall the proof here.  

\begin{proof}
First note that any $\Sigma $-comodule $N$ is the kernel of a map of
extended comodules.  Indeed, the structure map $\psi_{N} \mathcolon
N\xrightarrow{}\Sigma \otimes _{B}N$ is a comodule map if we give
$\Sigma \otimes _{B}N$ the extended comodule structure, and an
embedding because it is split over $B$ by the counit of $\Sigma $.  If
we let $q\mathcolon \Sigma \otimes _{B}N\xrightarrow{}Q$ denote the
quotient, then we get a diagram of comodules
\[
N \xrightarrow{\psi_{N} } \Sigma \otimes _{B} N \xrightarrow{\psi _{Q}q}
\Sigma \otimes _{B}Q
\]
expressing $N$ as the kernel of a map of extended comodules.
Adjointness forces us to define $\Phi^{*}(\Sigma \otimes _{B}P)=\Gamma \otimes
_{A}P$ for any $B$-module $P$.  Once we define $\Phi^{*}$ on maps between
extended comodules such as $\psi _{Q}q$, we can then define $\Phi^{*}N$ as the
kernel of $\Phi^{*}(\psi _{Q}q)$.  

So suppose we have a map $f\mathcolon \Sigma \otimes
_{B}P\xrightarrow{}\Sigma \otimes _{B}P'$.  We need to define
\[
\Phi^{*}f\mathcolon \Gamma \otimes _{A}P\xrightarrow{}\Gamma \otimes _{A}P'.
\]
By adjointness, it suffices to define a map of $A$-modules $\Gamma
\otimes _{A}P \xrightarrow{}P'$.  We define this map as the composite
\[
\Gamma \otimes _{A}P \xrightarrow{}\Sigma \otimes _{B}P
\xrightarrow{f}\Sigma \otimes _{B}P' \xrightarrow{}P'. 
\]
Here the first map is induced by the map $\Gamma \xrightarrow{}\Sigma $
and the last map is induced by the counit of $\Sigma $.  
\end{proof}
\begin{remark}
 It can be shown that when $N$ is a $\Sigma$-comodule, the group
$N\otimes_A\Gamma$ has compatible structures as a $\Gamma$-comodule and
a $\Sigma$-comodule.  This makes the $\Sigma$-primitives
$\Prim_\Sigma(N\otimes_A\Gamma)$ into a $\Gamma$-comodule, which turns
out to be isomorphic to $\Phi^*N$.  One can give a proof of the
existence of $\Phi^*$ based on this formula, but we do not need it so we
omit the details.
\end{remark}
We can now prove the main result of this section, which is also
Theorem~\ref{main-A} of the introduction.  

\begin{theorem}\label{thm-Giraud}
Suppose $(A, \Gamma )$ is a flat Hopf algebroid, and $B$ is a Landweber
exact $A$-algebra.  Let $\Phi \mathcolon
(A,\Gamma)\xrightarrow{}(B,\Gamma_{B})$ denote the corresponding map of
Hopf algebroids, inducing $\Phi_{*}\mathcolon \Gamma \comod
\xrightarrow{}\Gamma_{B}\comod$ with right adjoint $\Phi^{*}$.  Then the
counit of the adjunction $\epsilon \mathcolon
\Phi_{*}\Phi^{*}M\xrightarrow{}M$ is a natural isomorphism.  Hence
$\Phi^{*}\Phi_{*}$ is localization with respect to the hereditary
torsion theory $\cat{T}=\ker \Phi_{*}$, and $\Phi^{*}$ defines an
equivalence of categories from $\Gamma_{B}\comod$ to $L_{\cat{T}}(\Gamma
\comod )$.
\end{theorem}

\begin{proof}
Since $B$ is Landweber exact, $\Phi_{*}$ is exact, so $\epsilon$ is a
natural transformation of left exact functors.  Since every comodule is
a kernel of a map of extended comodules, it suffices to check that
$\epsilon_{N}$ is an isomorphism for extended comodules $N=\Gamma
_{B}\otimes_{B}V$.  But then we have
\begin{gather*}
\Phi_{*}\Phi^{*}N\cong B\otimes_{A} (\Gamma \otimes_{A} V)\cong
B\otimes_{A} \Gamma \otimes_{A}B\otimes_{B} V \\
\cong \Gamma _{B} \otimes_{B}V\cong M,
\end{gather*}
as required.  Proposition~\ref{prop-local-torsion} completes the proof. 
\end{proof}

\section{Torsion theories of $BP_{*}BP$-comodules}
 
Suppose $(A,\Gamma )$ is a flat Hopf algebroid, and $B$ is a Landweber
exact $A$-algebra.  Theorem~\ref{thm-Giraud} shows that the category of
$\Gamma _{B}$-comodules is equivalent to the localization of the
category of $\Gamma $-comodules with respect to some hereditary torsion
theory $\cat{T}$.  Thus we would like to classify all hereditary torsion
theories of $\Gamma $-comodules.  This is of course impossible in
general, but it turns out to be tractable in the cases of most interest
in algebraic topology.  In this section, we concentrate on the case $(A,
\Gamma)=(BP_{*},BP_{*}BP)$, where $BP$ is the Brown-Peterson spectrum.

Recall that $BP_{*}\cong \Zp [v_{1},v_{2},\dotsc]$, where
$|v_{i}|=2(p^{i}-1)$.  We choose the $v_{i}$ to be the Araki
generators~\cite[Section~A2.2]{ravenel} for definiteness, but all that
matters is that $v_{n}$ is primitive modulo $I_{n}=(p,v_{1},\dotsc
,v_{n-1})$.  The ideals $I_{n}$ are independent of the choice of
generators.  For notational purposes, we take $v_{0}=p$ and $v_{-1}=0$.  
We also write $s$ for the shift functor on $BP_{*}BP$-comodules, so that
$(sM)_{n}=M_{n-1}$.  

Let $\cat{T}_{n}$ denote the class of all graded $BP_{*}BP$-comodules
that are $v_{n}$-torsion, in the sense that each element is killed by
some power of $v_{n}$, depending on the element.  By Lemma~2.3
of~\cite{johnson-yosimura}, $M$ is $v_{n}$-torsion if and only if $M$ is
$I_{n+1}$-torsion, so $\cat{T}_{n}$ is independent of the choice of
generators.  

The following theorem is Theorem~\ref{main-B} of the introduction.  

\begin{theorem}\label{thm-torsion}
Let $\cat{T}$ be a graded hereditary torsion theory of graded
$BP_{*}BP$-comodules, and suppose that $\cat{T}$ contains some nonzero
comodule that is finitely presented over $BP_{*}$.  Then
$\cat{T}=\cat{T}_{n}$ for some $n\geq -1$.  
\end{theorem}

The reader should compare Theorem~\ref{thm-torsion} to the
classification of Serre classes of finitely presented
$BP_{*}BP$-comodules in~\cite{jeanneret-landweber-ravenel} (which they
call thick subcategories).  Given a hereditary torsion theory
$\cat{T}$, the collection $\cat{F}$ of all finitely presented comodules
in it is a Serre class (of all the finitely presented comodules);
combining Theorem~\ref{thm-torsion} with the result
of~\cite{jeanneret-landweber-ravenel} says that as long as $\cat{F}$ is
nonzero, then $\cat{T}$ is uniquely determined by $\cat{F}$.

We do not know what happens when there are no nonzero finitely presented
comodules in $\cat{T}$.  In this case, Proposition~\ref{prop-non-nil}
below implies that every comodule in $\cat{T}$ is $v_{n}$-torsion for
all $n$.  Ravenel~\cite[Section~2]{rav-loc} conjectures that there are
uncountably many different Bousfield classes of spectra $BPI$ where $I$
is an infinite regular sequence in $BP_{*}$.  One might similarly
conjecture that there are uncountably many different hereditary torsion
theories $\cat{T}$ containing no nonzero finitely presented comodules.

Theorem~\ref{thm-torsion} will follow from the two propositions below.  

\begin{proposition}\label{prop-nil}
$\cat{T}_{n}$ is the graded hereditary torsion theory generated by the
$BP_{*}BP$-comodule $BP_{*}/I_{n+1}$.
\end{proposition}

\begin{proposition}\label{prop-non-nil}
Suppose that $\cat{T}$ is a graded hereditary torsion theory of graded
$BP_{*}BP$-comodules such that $BP_{*}/I_{n}\not \in \cat{T}$.  Then 
$\cat{T}\subseteq \cat{T}_{n}$.
\end{proposition}

Given these two propositions, Theorem~\ref{thm-torsion} follows easily.  

\begin{proof}[Proof of Theorem~\ref{thm-torsion}] 
Suppose $\cat{T}$ is a graded hereditary torsion theory containing the
nonzero finitely presented comodule $M$.  The Landweber filtration
theorem for $BP_{*}BP$-comodules~\cite[Theorem~2.3]{land-exact}
guarantees that $M$ has a subcomodule of the form $s^{t}BP_{*}/I_{j}$
for some $j$ and some $t$.  Thus $BP_{*}/I_{j}\in \cat{T}$.  Let
\[
n+1=\min \{j\suchthat BP_{*}/I_{j}\in \cat{T} \}.
\]
Then $\cat{T}\supseteq \cat{T}_{n}$ by Proposition~\ref{prop-nil}.  On
the other hand, $BP_{*}/I_{n}\not \in \cat{T}$, so $\cat{T}\subseteq
\cat{T}_{n}$ by Proposition~\ref{prop-non-nil}.  Hence
$\cat{T}=\cat{T}_{n}$, as required.
\end{proof}

We owe the reader proofs of Proposition~\ref{prop-nil} and
Proposition~\ref{prop-non-nil}.  We need the following lemma.

\begin{lemma}\label{lem-nil}
Suppose $M$ is a nonzero $v_{n}$-torsion graded $BP_{*}BP$-comodule.
Then $M$ has a nonzero primitive $x$ such that $I_{n+1}\subseteq \Ann
(x)$.
\end{lemma}

\begin{proof}
Let $y$ be a nonzero element of $M$, and let $I=\sqrt{\Ann y}$.  Since
$y$ is $v_{n}$-torsion, it is also $I_{n+1}$-torsion, and so
$I_{n+1}\subseteq I$.  Theorem~2 of~\cite{landweber-comodules}
guarantees that there is a primitive $x$ with $\Ann (x)=I$. 
\end{proof}

\begin{proof}[Proof of Proposition~\ref{prop-nil}]
Let $\cat{T}$ denote the graded hereditary torsion theory generated by
$BP_{*}/I_{n+1}$.  Since one can easily check that $\cat{T}_{n}$ is a
graded hereditary torsion theory, and $BP_{*}/I_{n+1}\in \cat{T}_{n}$,
we see that $\cat{T}\subseteq \cat{T}_{n}$.  Conversely, suppose $M$ is
$v_{n}$-torsion.  We construct a transfinite increasing sequence
$M_{\alpha}$ of subcomodules of $M$ such that each $M_{\alpha}$ is in
$\cat{T}$.  This sequence will be strictly increasing unless
$M_{\beta }=M$ for some $\beta $, so we conclude that $M=M_{\beta }\in
\cat{T}$.

To construct $M_{0}$, we use Lemma~\ref{lem-nil} to find a nonzero
primitive $x\in M$ such that $I_{n+1}x=0$.  This gives a
subcomodule $M_{0}\cong s^{t}BP_{*}/I$ of $M$ such that $I\supseteq
I_{n+1}$.  Hence $M_{0}$ is isomorphic to a quotient of
$s^{t}BP_{*}/I_{n+1}$, so $M_{0}\in \cat{T}$.  This begins the
transfinite induction.  The limit ordinal step of the induction is
simple.  If we have defined $M_{\alpha}$ for all $\alpha <\beta$ for
some limit ordinal $\beta$, we define $M_{\beta}=\bigcup_{\alpha
<\beta}M_{\alpha}$.  Since $M_{\beta}$ is a quotient of
$\bigoplus_{\alpha <\beta}M_{\alpha}$, $M_{\beta}$ is still in
$\cat{T}$.

We now carry out the successor ordinal step of the induction.  So
suppose we have defined $M_{\alpha}$.  If $M_{\alpha}=M$, we let
$M_{\alpha +1}=M$ as well.  Otherwise, consider the quotient
$M/M_{\alpha}$.  Since this comodule is $v_{n}$-torsion,
Lemma~\ref{lem-nil} gives us an element $z\in M$ such that the coset
$\overline{z}$ in $M/M_{\alpha}$ is a nonzero primitive such that
$I_{n+1}\overline{z}=0$.  We define $M_{\alpha+1}$ to be the subcomodule
of $M$ generated by $M_{\alpha}$ and $z$.  Then $M_{\alpha +1}$ is an
extension of $M_{\alpha}$ and $s^{t}BP_{*}/\Ann (\overline{z})$ for
some $t$, so $M_{\alpha +1}\in \cat{T}$ as required.  This completes the
proof.
\end{proof}

\begin{proof}[Proof of Proposition~\ref{prop-non-nil}]
Suppose $BP_{*}/I_{n}$ is not in some graded hereditary torsion theory
$\cat{T}$, and $M$ is in $\cat{T}$.  We must show that every $x\in M$ is
$v_{n}$-torsion.  We show that every $x\in M$ is $v_{j}$-torsion for
$-1\leq j\leq n$ by induction on $j$.  The initial step is automatic
since $v_{-1}=0$.  For the induction step, assume that we have shown
that $M$ is $v_{j-1}$-torsion, and let $M_{j}$ denote the
$v_{j}$-torsion in $M$.  This is a subcomodule of $M$ by Proposition~2.9
of~\cite{johnson-yosimura}.  Suppose that $M_{j}$ is not all of $M$.
Lemma~\ref{lem-nil} allows us to find a primitive $y_{j}$ in $M/M_{j}$
with $I_{j}\subseteq \Ann (y_{j})$. Since $M/M_{j}$ is
$v_{j}$-torsion-free, and the only primitives modulo $I_{j}$ are powers
of $v_{j}$ (\cite[Theorem~4.3.2]{ravenel}), we conclude that $\Ann
(y_{j})=I_{j}$.  This gives us an embedding $s^{t}BP_{*}/I_{j}\subseteq
M/M_{j}$, so $BP_{*}/I_{j}\in \cat{T}$.  This contradiction shows that
$M$ is all $v_{j}$-torsion, as required.
\end{proof}

\section{Equivalences of comodule categories}\label{sec-results}

In this section, we show that the category of comodules over a
Landweber exact $BP_{*}$-algebra $B$ depends only on the height of
$B$, and deduce versions of the Miller-Ravenel and Morava change of
rings theorems.

\begin{definition}\label{defn-height}
Suppose $B$ is a nonzero graded $BP_{*}$-module.  We define the
\textbf{height} of $B$, written $\height B$, to be the largest $n$ such that
$B/I_{n}$ is nonzero, or $\infty$ if $B/I_{n}$ is nonzero for all $n$.
\end{definition}

Note that $\height E(n)_{*}=\height v_{n}^{-1}BP_{*}=n$ and $\height
BP_{*}=\infty$.  The following theorem implies Theorem~\ref{main-C} of
the introduction.  

\begin{theorem}\label{thm-main-BP}
Let $(A, \Gamma )=(BP_{*}, BP_{*}BP)$, and suppose $B$ and $B'$ are two
graded Landweber exact $BP_{*}$-algebras with $\height B=\height B'=n
\leq \infty $.  Then the category of graded $\Gamma _{B}$-comodules is
equivalent to the category of graded $\Gamma _{B'}$-comodules.  If
$n=\infty $, these categories are equivalent to the category of graded
$\Gamma $-comodules.  If $n<\infty $, these categories are equivalent to
the localization of the category of graded $\Gamma $-comodules with
respect to the graded hereditary torsion theory $\cat{T}_{n}$.
\end{theorem}

\begin{proof}
Theorem~\ref{thm-Giraud} implies that the category of graded $\Gamma
_{B}$-comodules is equivalent to the localization of the category of
graded $\Gamma $-comodules with respect to the kernel $\cat{T}$ of the
functor $M\mapsto B\otimes _{A}M$.  Suppose first that $n<\infty $.
Then $BP_*/I_{n}\not \in \cat{T}$ but $BP_{*}/I_{n+1}\in \cat{T}$.
Theorem~\ref{thm-torsion}, or, more precisely,
Propositions~\ref{prop-nil} and~\ref{prop-non-nil}, imply that
$\cat{T}=\cat{T}_{n}$.  

Now suppose that $n=\infty $.  We claim that $\cat{T}=(0)$.  Indeed,
suppose $M\in \cat{T}$ and is nonzero.  Since every graded
$BP_{*}BP$-comodule has a primitive, $M$ has a subcomodule isomorphic
to $s^{t}BP_{*}/I$ for some invariant ideal $I$.  But $BP_{*}$ is a
local ring (in the graded sense), with maximal ideal
$I_{\infty}=\bigcup_{r} I_{r}$.  Thus $I\subseteq I_{\infty}$, and so
$BP_{*}/I_{\infty}\in \cat{T}$, since it is a quotient of $BP_{*}/I$.
Hence $B/(I_{\infty}B)=0$, and so the unit $1\in B$ is in
$I_{\infty}B$, so must be in $I_{r}B$ for some $r$.  Hence
$B/I_{r}B=0$, contradicting our assumption that $B$ has infinite
height.  Thus $\cat{T}=(0)$.
\end{proof}

In view of this theorem, we denote the localization functor on
$BP_{*}BP$-comodules corresponding to the torsion theory $\cat{T}_{n}$
by $L_{n}\mathcolon BP_{*}BP\comod \xrightarrow{} BP_{*}BP\comod$.  

\begin{corollary}\label{cor-main-BP}
Let $(A, \Gamma )=(BP_{*}, BP_{*}BP)$, and suppose $B$ is a Landweber
exact $BP_{*}$-module of finite height $n$. 
Then the category of graded $\Gamma _{B}$-comodules is equivalent to the
full subcategory of graded $\Gamma $-comodules $M$ such that 
\[
\Hom _{A}^{*}(A/I_{n+1},M) = \Ext ^{1,*}_{\Gamma } (A/I_{n+1},M) =0
\]
\end{corollary}

\begin{proof}
It suffices to prove that $M$ is $L_{n}$-local if and only if $M$
satisfies the given conditions.  Lemma~\ref{lem-E-local} tells
us that $M$ is local if and only if
\[
\Hom _{\Gamma }^{*}(T,M) = \Ext ^{1,*}_{\Gamma }(T,M)=0
\]
for all $v_{n}$-torsion comodules $T$.  

We first claim that $\Hom _{\Gamma }^{*}(T,M)=0$ for all $v_{n}$-torsion
comodules $T$ if and only if $\Hom _{A}^{*}(A/I_{n+1},M)=0$.  To see
this, note that $\Hom _{A}^{*}(A/I_{n+1},M)=0$ if and only if $M$ has no
nonzero $v_{n}$-torsion elements (using the equivalence of
$v_{n}$-torsion with
$I_{n+1}$-torsion~\cite[Lemma~2.3]{johnson-yosimura}).  But then $\Hom
_{A}^{*}(T,M)=0$ for all $v_{n}$-torsion modules $T$, and so, a fortiori,
$\Hom _{\Gamma }^{*}(T,M)=0$ for all $v_{n}$-torsion comodules $T$.
Conversely, if $\Hom _{\Gamma }^{*}(A/I_{n+1},M)=0$, then $M$ has no
primitives that are $v_{n}$-torsion.  Proposition~2.7
of~\cite{johnson-yosimura} implies that $M$ has no $v_{n}$-torsion at
all, and so $\Hom _{A}^{*}(A/I_{n+1},M)=0$ as well.

We now claim that $\Ext ^{1,*}_{\Gamma }(A/I_{n+1},M)=0$ if and only if
$\Ext ^{1,*}_{\Gamma }(T,M)=0$ for all $v_{n}$-torsion comodules $T$, as
long as $\Hom _{A}^{*}(A/I_{n+1},M)=0$.  Indeed, suppose $\Ext
^{1,*}_{\Gamma }(A/I_{n+1},M)=0$.  Using the short exact sequences
\[
0 \xrightarrow{} A/I_{n+j} \xrightarrow{v_{n+j}} A/I_{n+j}
\xrightarrow{} A/I_{n+j+1} \xrightarrow{} 0,
\]
we conclude that $\Ext ^{1,*}_{\Gamma }(A/I_{r},M)=0$ for all $r>n$.
The Landweber filtration theorem~\cite[Theorem~2.3]{land-exact} tells us
that every finitely presented $\Gamma $-comodule $T$ is an iterated
extension of suspensions of comodules of the form $A/I_{r}$; if $T$ is
$v_{n}$-torsion, one can easily check that $r>n$.  Thus $\Ext
^{1}_{\Gamma }(T,M)=0$ for all finitely presented $v_{n}$-torsion
comodules.  Now every comodule $T$ is a filtered colimit of finitely
presented comodules, by Lemma~1.15 of~\cite{johnson-yosimura}.  If $T$
is $v_{n}$-torsion, then the finitely presented $v_{n}$-torsion
comodules mapping to $T$ are cofinal, so $T$ is a filtered colimit
$\colim T_{\alpha }$ of finitely presented $v_{n}$-torsion comodules.
This gives a short exact sequence
$T''\xrightarrow{}T'\xrightarrow{}T$, where 
$T'=\bigoplus_\alpha T_\alpha$ (so $\Ext^{1,*}_\Gamma(T',M)=0$) and
$T''\leq T'$ (so $T''$ is $v_n$-torsion and thus
$\Hom^*_{\Gamma}(T'',M)=0$).  It follows that
$\Ext^{1,*}_\Gamma(T,M)=0$ as required.

\end{proof}

In~\cite{hovey-strickland-derived}, we strengthen
Corollary~\ref{cor-main-BP} by showing that a $BP_{*}BP$-comodule $M$ is
$L_{n}$-local if and only if 
\[
\Hom_{BP_{*}}^{*}(BP_{*}/I_{n+1},M) =
\Ext^{1,*}_{BP_{*}}(BP_{*}/I_{n+1},M)=0,
\]
with the $\Ext$ group computed over the ring $BP_*$ rather than the
Hopf algebroid $BP_*BP$.  This $\Ext$ group can easily be computed
from a Koszul resolution and so is much more accessible than the
previous one.

\begin{corollary}\label{cor-maps}
Let $(A,\Gamma )=(BP_{*},BP_{*}BP)$, and suppose $B\xrightarrow{}B'$ is
a map of Landweber exact $A$-algebras.  Let $\cat{T}$ denote the graded
hereditary torsion theory of $\Gamma _{B}$-comodules generated by
$B/I_{\height B'+1}$ if $\height B'<\infty $ and $(0)$ if
$\height B'=\infty $.  Then the functor $M\mapsto B'\otimes _{B}M$
defines an equivalence of categories between the localization of the
category of graded $\Gamma _{B}$-comodules with respect to $\cat{T}$ and
the category of graded $\Gamma _{B'}$-comodules.  In particular, if
$\height B=\height B'$, then this functor is itself an equivalence of
categories.
\end{corollary}

This corollary is a special case of the general fact that maps between
localizations are themselves localizations;
see~\cite[Lemma~3.1.5]{hovey-axiomatic}.  

\begin{example}\label{ex-exact}
There are well-known maps 
\[
v_{n}^{-1}BP_{*} \xrightarrow{} E(n)_{*} \xrightarrow{} E_{n*}
\]
of Landweber exact $BP_{*}$-algebras of height $n$.  These
maps induce equivalences of the associated categories of comodules.
Note that they certainly do not induce equivalences of the associated
categories of modules; in particular, $E(n)_{*}$ is Noetherian and
$v_{n}^{-1}BP_{*}$ is not.
\end{example}

We can now give a straightforward and systematic account of some
well-known change of rings theorems, as mentioned in the introduction.
The following is our general result; it follows immediately from
Corollary~\ref{cor-maps}.

\begin{theorem}\label{thm-change-of-rings}
Let $(A, \Gamma )=(BP_{*},BP_{*}BP)$, and suppose $B\xrightarrow{}B'$ is
a map of Landweber exact $A$-algebras such that $\height B=\height B'$.
Then
\[
\Ext _{\Gamma _{B}}^{**}(M,N) \cong \Ext _{\Gamma
_{B'}}^{**}(B'\otimes _{B}M, B'\otimes _{B}N).   \qed
\]
\end{theorem}

The Morava change of rings theorem~\cite{morava-comodules} is often
stated in precisely this form.  We give a graded version of it, as
opposed to the ungraded version given in~\cite[Theorem~6.1.3]{ravenel}.  

\begin{corollary}\label{cor-morava}
Suppose $(A, \Gamma )=(BP_{*},BP_{*}BP)$, and let $I$ denote the
ideal in $A$ generated by $p$ and all the $v_{i}$ except
$v_{n}$.  Let $B$ denote the completion of $v_{n}^{-1}A$ at $I$, and let
$B'$ denote the completion of $E(n)_{*}$ at $I_{n}$.  Then
\[
\Ext _{\Gamma _{B}}^{**}(M,N) \cong \Ext _{\Gamma _{B'}}^{**}(B'\otimes
_{B}M, B'\otimes _{B}N) 
\]
for all $\Gamma _{B}$-comodules $M$ and $N$.  
\end{corollary}

Note that the Morava change of rings theorem was only known before in
case $M=B$.  

Here is our version of the Miller-Ravenel change of rings
theorem~\cite[Theorem~3.10]{miller-ravenel}.

\begin{corollary}\label{cor-miller-ravenel}
Let $(A, \Gamma )=(BP_{*},BP_{*}BP)$, $B=v_{n}^{-1}BP_{*}$, and
$B'=E(n)_{*}$.  Then
\[
\Ext _{\Gamma _{B}}^{**}(M,N) \cong \Ext _{\Gamma _{B'}}^{**}(B'\otimes
_{B}M, B'\otimes _{B}N)
\]
for all $\Gamma _{B}$-comodules $M$ and $N$.  
\end{corollary}

The Miller-Ravenel change of rings theorem is usually stated as
\begin{equation}\label{eq-miller-ravenel}
\Ext _{BP_{*}BP}^{**}(BP_{*},N) \cong  \Ext
_{E(n)_{*}E(n)}^{**}(E(n)_{*}, E(n)_{*}\otimes _{BP_{*}}N)
\end{equation}
for all $BP_{*}$-comodules $N$ on which $v_{n}$ acts invertibly.  This
is a consequence of Corollary~\ref{cor-miller-ravenel}, arguing as in
Lemmas~3.11 and~3.12 of~\cite{miller-ravenel}.  The point is
essentially as follows: if $v_n$ acts invertibly on $N$, then nothing
changes if we invert $v_n$ in $BP_*$ and $BP_*BP$.  Moreover, $N$ is
necessarily $I_n$-torsion, so $v_n$ is asymptotically primitive on
anything involving $N$, so we need not distinguish between inverting
$\eta_L(v_n)$ or $\eta_R(v_n)$.  

We have also generalized the statement~\ref{eq-miller-ravenel} of the
Miller-Ravenel change of rings theorem in~\cite{hovey-barcelona},
expressing it as an isomorphism between morphism sets in appropriate
derived categories.

Similarly, if $m\geq n$, we can apply Theorem~\ref{thm-change-of-rings}
to the map $v_{n}^{-1}BP_{*}\xrightarrow{}v_{n}^{-1}E(m)_{*}$ to get a
version of the change of rings theorem
of~\cite[Theorem~3.1]{hovey-sadofsky-picard}.

\begin{corollary}\label{cor-hovey-sadofsky}
Suppose $(A, \Gamma )=(BP_{*},BP_{*}BP)$, and let
$B=v_{n}^{-1}BP_{*}$ and $B'=v_{n}^{-1}E(m)_{*}$ for $m\geq n$.  Then 
\[
\Ext _{\Gamma _{B}}^{**}(M,N) \cong \Ext _{\Gamma _{B'}}^{**}(B'\otimes
_{B}M, B'\otimes _{B}N)
\]
for all $\Gamma _{B}$-comodules $M$ and $N$.  
\end{corollary}

Again, the methods of Miller and Ravenel allow one to derive the
original change of rings theorem of the first author and Sadofsky from
Corollary~\ref{cor-hovey-sadofsky}.

\section{The structure of $E_{*}E$-comodules}\label{sec-structure}

This section is devoted to proving analogues for $E_{*}E$-comodules of
the Landweber structure theorems for $BP_{*}BP$-comodules, when $E_{*}$ is
Landweber exact over $BP_{*}$.

\begin{theorem}\label{thm-primitive}
 Let $(A, \Gamma )=(BP_{*}, BP_{*}BP)$, and suppose $B$ is a Landweber
 exact $A$-algebra.  Then every nonzero $\Gamma _{B}$-comodule has a
 nonzero primitive.
\end{theorem}

\begin{proof}
Let $\Phi \mathcolon (A, \Gamma)\xrightarrow{}(B, \Gamma_{B})$ denote
the evident map of Hopf algebroids.  Let $\Phi_{*}\mathcolon \Gamma
\comod \xrightarrow{}\Gamma_{B}\comod$ denote the induced functor, with
right adjoint $\Phi^{*}$.  Suppose $M$ is a nonzero
$\Gamma_{B}$-comodule.  We must show that $\Hom_{\Gamma_{B}}^{*}(B,M)$
is nonzero.  But adjointness implies that 
\[
\Hom_{\Gamma_{B}}^{*}(B,M)=\Hom_{\Gamma_{B}}^{*}(\Phi_{*}A,M)\cong
\Hom_{\Gamma}^{*}(A, \Phi^{*}M).  
\]
Since $\Phi_{*}\Phi^{*}M\cong M$ by Theorem~\ref{thm-Giraud},
$\Phi^{*}M$ is a nonzero $\Gamma$-comodule.  It is well-known that every nonzero
$\Gamma$-comodule has a primitive; it follows for example from
Lemma~\ref{lem-nil}.  Thus $\Hom_{\Gamma}^{*}(A, \Phi^{*}M)$ is nonzero,
as required.  
\end{proof}

We now compute the primitives in $B/I_{n}$ for all $n$.  The case
$B=BP_{*}$ is well-known~\cite[Theorem~4.3.2]{ravenel}.  

\begin{theorem}\label{thm-calculation}
Let $(A, \Gamma )=(BP_{*},BP_{*}BP)$, and suppose $B$ is a nonzero
Landweber exact $A$-algebra.
\begin{enumerate}
\item [(a)] If $\height B>0$, then $\Hom^{*}_{\Gamma_{B}}(B,B)\cong
\Zp$.  
\item [(b)] If $\height B=0$, then $\Hom^{*}_{\Gamma_{B}}(B,B)\cong
\Q$.  
\item [(c)] If $\height B > n>0$, then $\Hom^{*}_{\Gamma_{B}}(B,
B/I_{n})\cong \Fp [v_{n}]$.  
\item [(d)] If $\infty>\height B=n>0$, then $\Hom^{*}_{\Gamma_{B}}(B,
B/I_{n})\cong \Fp [v_{n}, v_{n}^{-1}]$.
\item [(e)] If $\height B=\infty$ then  $\Hom^{*}_{\Gamma_{B}}(B,
B/I_{\infty})\cong \Fp$.
\item [(f)] If $n>\height B$ then $B/I_n=0$ and so
 $\Hom^{*}_{\Gamma_{B}}(B,B/I_n)=0$.
\end{enumerate}
\end{theorem}

The simplest way to prove this theorem is to use the following
computation. Recall that $L_{n}$ denotes the localization functor on the
category of $BP_{*}BP$-comodules with respect to the hereditary torsion
theory of $v_{n}$-torsion comodules.  

\begin{lemma}\label{lem-calculation}
For $n<\infty$ we have 
\[
L_{n}(BP_{*}/I_{n})\cong v_{n}^{-1}BP_{*}/I_{n} 
\]
and 
\[
L_{n}(BP_{*}/I_{m}) = BP_{*}/I_{m}
\]
for $m<n$. 
\end{lemma}

As usual, we let $v_{0}=p$ and $I_{0}=(0)$ in interpreting the statement
of this lemma.  Recall also that $L_\infty$ is the identity functor,
so the $n=\infty$ case is trivial.

\begin{proof}
Let $M$ denote either $BP_{*}/I_{m}$ (for $m<n$) or
$v_{n}^{-1}BP_{*}/I_{n}$ (for $m=n$).  It suffices to show that $M$ is
$L_{n}$-local, since the map
$BP_{*}/I_{n}\xrightarrow{}v_{n}^{-1}BP_{*}/I_{n}$ has $v_{n}$-torsion
cokernel, so is an $L_{n}$-equivalence.  For this, we use
Corollary~\ref{cor-main-BP}.  Since $M$ has no $v_{n}$-torsion, it
suffices to show that $\Ext^{1,*}_{BP_{*}BP}(BP_{*}/I_{n+1}, M)=0$.  

We first show that $\Ext^{1,*}_{BP_{*}}(BP_{*}/I_{n+1},M)=0$.  So
suppose we have a short exact sequence of $BP_{*}$-modules
\begin{equation}\label{eq-short-exact}
0 \xrightarrow{} M \xrightarrow{i} X \xrightarrow{g} s^{t}BP_{*}/I_{n+1}
\xrightarrow{} 0.  
\end{equation}
Let $x\in X$ be such that $g(x)$ is the generator of
$s^{t}BP_{*}/I_{n+1}$.  The argument now depends on whether
$M=BP_{*}/I_{m}$ or $M=v_{n}^{-1}BP_{*}/I_{n}$.   In case
$M=v_{n}^{-1}BP_{*}/I_{n}$, note that $v_{n}x\in M$.  Since $v_{n}$ acts
invertibly on $M$, there is a $y\in M$ such that $v_{n}y=v_{n}x$.  Then
$v_{n}(x-y)=0$, and also $v_{i}(x-y)\in M$ for $i<n$.  Hence
$v_{n}v_{i}(x-y)=0$, and so $v_{i}(x-y)=0$ since $M$ has no
$v_{n}$-torsion.  Thus $x-y$ defines a splitting of our given exact
sequence~\ref{eq-short-exact}.  

Now suppose $M=BP_{*}/I_{m}$ for some $m<n$.  Then $v_{m}x$ and
$v_{m+1}x$ are both in $M$.  Since $v_{m+1}(v_{m}x)=v_{m}(v_{m+1}x)$ and
$M$ is a unique factorization domain, we conclude that $v_{m}x=v_{m}y$
for some $y\in M$, and that $v_{m+1}x=v_{m+1}y$.  Now a similar argument
as we used in case $M=v_{n}^{-1}BP_{*}/I_{n}$ shows that $v_{i}x=v_{i}y$
for all $i\leq n$.  Hence $x-y$ defines a splitting
of~\ref{eq-short-exact}.  

Now suppose the short exact sequence~\ref{eq-short-exact} is a sequence
of $BP_{*}BP$-comodules.  Write $X\cong M\oplus s^{t}BP_{*}/I_{n+1}$ as
$BP_{*}$-modules.  We claim that this splitting must be a splitting of
$BP_{*}BP$-comodules as well.  Indeed, let $Y$ be the $v_{n}$-torsion
in $X$, which is a subcomodule and is obviously just 
$0\oplus s^{t}BP_{*}/I_{n+1}$.  Hence the map 
$Y\xrightarrow{} X\xrightarrow{}BP_*/I_{n+1}$ is an isomorphism, and
its inverse gives a splitting of the sequence.
\end{proof}

\begin{proof}[Proof of Theorem~\ref{thm-calculation}]
As usual, let $\Phi_{*}$ denote the functor from $\Gamma$-comodules to
$\Gamma_{B}$-comodules that takes $M$ to $B\otimes_{A}M$.  Then we have 
\begin{gather*}
\Hom_{\Gamma_{B}}^{*}(B, B/I_{m}) =\Hom_{\Gamma_{B}}^{*}(\Phi_{*}A,
\Phi_{*}A/I_{m}) \\
\cong \Hom_{\Gamma}^{*}(A, \Phi^{*}\Phi_{*}(A/I_{m})) \cong
\Hom_{\Gamma}^{*}(A,L_{n}(A/I_{m})).
\end{gather*}
Parts~(a) to~(d) of the theorem now follow from
Lemma~\ref{lem-calculation} and~\cite[Proposition~5.1.12]{ravenel}.
Part~(e) follows from the observation that
$\Hom_\Gamma^*(A,BP_*/I_\infty)=BP_*/I_\infty=\Fp$.
Part~(f) is just the definition of $\height B$, recorded for ease of
comparison. 
\end{proof}

We now consider the analogue of Landweber's classification of invariant
radical ideals in $BP_*$~\cite[Theorem~2.2]{land-exact}.  For this to
work smoothly, we need to modify the problem slightly.  Consider an
abelian category $\cat{A}$ with a symmetric monoidal tensor product, and
let $k$ be the unit for the tensor product.  We define a
\textbf{categorical ideal} in $\cat{A}$ to be a subobject of $k$; given
ideals $I$ and $J$, we let $IJ$ denote the image of the evident map
$I\otimes J\xrightarrow{}k$.  We say that $I$ is \textbf{categorically
radical} if $J^2\leq I$ implies $J\leq I$.  This notion is evidently
invariant under monoidal equivalences of abelian categories, such as
those in Theorem~\ref{thm-main-BP}.

We now specialize to the case $\cat{A}=(B,\Sigma)\comod$.  The
categorical ideals are then the invariant ideals in $B$.  An invariant
radical ideal is categorically radical, but the converse need not be
true.  For example, take $(A,\Gamma)=(BP_*,BP_*BP)$ as before, and
$B=BP_*[x]/x^2$, and $\Sigma=\Gamma_B$.  Then $I_nB$ is categorically
radical, but not radical.  Indeed, the only invariant ideals are those
of the form $IB$ for some invariant ideal $I\leq BP_*$, and $IB$ is
never radical.

One can easily check that the proof of Landweber's classification of
invariant radical ideals in $BP_{*}$ in~\cite[Theorem~4.3.2]{ravenel}
also classifies categorically radical ideals.  That is, we have the
following theorem.  

\begin{theorem}\label{thm-categorically-invariant}
The ideal $I\leq BP_{*}$ is a categorically radical invariant ideal if
and only if $I=I_{n}$ for some $0\leq n\leq \infty $.
\end{theorem}

Almost the same theorem holds for categorically radical ideals in
Landweber exact $BP_{*}$-algebras.  

\begin{theorem}\label{thm-invariant}
 Suppose $(A, \Gamma )=(BP_{*},BP_{*}BP)$ and $B$ is a Landweber exact
 $A$-algebra.  Then the categorically radical invariant ideals in $B$
 are $\{I_kB\suchthat 0\leq k\leq\height B\}$.  In particular, this set
 contains all the invariant radical ideals.
\end{theorem}

\begin{proof}
 Put $n=\height B$.  We assume that $n>0$, leaving the rational case to
the reader.  As $n>0$, we have $\Phi^*B=L_nBP_*=BP_*=A$.  Note also that
$\Phi^*$ is left exact, so it preserves monomorphisms, so it sends
invarant ideals in $B$ to invariant ideals in $A$.  Consider a
categorically radical invariant ideal $J\leq B$, and put $K=\Phi^*J\leq
A$, so $J=\Phi_*K=KB$.  If $K=A$ this means that $J=B$, which we have
implicitly excluded from consideration; so $K$ is a proper ideal in $A$.
We claim that $K$ is also categorically radical.  Indeed, suppose we
have an invariant ideal $K_0$ with $K_0^2\leq K$.  Put $J_0=BK_0$ and
apply $\Phi_*$ to the maps $K_0\otimes
K_0\xrightarrow{}K\xrightarrow{}A$ to see that $J_0^2\leq J$.  As $J$ is
categorically radical, we have $J_0\leq J$, so
$L_nK_0=\Phi^*\Phi_*K_0=\Phi^*J_0\leq\Phi^*J=K$.  Moreover, as $A$ is an
integral domain, $K_0$ has no $I_{n+1}$-torsion, so $K_0\leq L_nK_0$, so
$K_0\leq K$ as required.

 Since $K$ is categorically radical, we must have $K=I_{k}$ for some
$0\leq k\leq \infty$.  If $k>\height B$, then $J=\Phi_{*}K=B$.  Hence
$0\leq k\leq \height B$.  
\end{proof}




This classification leads to the analogue of the Landweber filtration
theorem, proved by Landweber~\cite{land-exact} for
$BP_{*}BP$-comodules.

\begin{theorem}\label{thm-filtration}
Suppose $(A, \Gamma )=(BP_{*},BP_{*}BP)$, and let $B$ be a Landweber
exact $A$-algebra.  Then every $\Gamma _{B}$-comodule $M$ that is
finitely presented over $B$ admits a finite filtration by subcomodules
\[
0 = M_{0} \subseteq M_{1} \subseteq \dots \subseteq M_{s} =M
\]
for some $s$, with $M_{r}/M_{r-1}\cong s^{t_{r}}B/I_{j_{r}}$ with
$j_{r}\leq \height B$ for all $r$.  
\end{theorem}

\begin{proof}
First note that $M$ is finitely presented over $B$ if and only if $M$ is
a finitely presented object of $\Gamma_{B}\comod$; that is, if and only
if $\Hom_{\Gamma_{B}}^{*}(M,-)$ commutes with all filtered colimits.
This is proved in~\cite[Proposition~1.3.3]{hovey-comodule-homotopy}.  It
follows that the statement of Theorem~\ref{thm-filtration} is invariant
under the equivalences of categories in Theorem~\ref{thm-main-BP}.
Hence we can assume that either $B=BP_{*}$ or $B=E(n)_{*}$.  The case
$B=BP_{*}$ is the classical Landweber filtration theorem.  

So now suppose $B=E(n)_{*}$ and $M$ is an arbitrary graded
$\Gamma_{B}$-comodule.  We construct a subcomodule of $M$ isomorphic to
$s^{t}B/I_{m}$ for some $p$ and some $m\leq n$.  Indeed, choose a
nonzero primitive $y_{0}$ in $M$.  If $\Ann (y_{0})=(0)$, we are done.
If not, $p^{i}y_{0}=0$ for some $i$ by Theorem~\ref{thm-calculation}.
Choose a minimal such $i$ and let $y_{1}=p^{i-1}y_{0}$.  Then $\Ann
(y_{1})$ is a proper invariant ideal containing $(p)$.  If it is $(p)$,
we are done.  If not, then Theorem~\ref{thm-calculation} implies that
$v_{1}^{j}y_{1}=0$ for some minimal $j>0$.  Let
$y_{2}=v_{1}^{j-1}y_{1}$.  Then $y_{2}$ is primitive (since $v_{1}$ is
primitive mod $p$), and $\Ann (y_{2})$ is an invariant ideal containing
$I_{2}$.  We continue in this fashion until we reach an $m$ such that
$\Ann (y_{m})=I_{m}$.  This must happen before we reach $n+1$.

Now we construct $M_{i}$ by induction, by applying this construction to
$M/M_{i-1}$.  Since $M$ is finitely generated and $B$ is Noetherian, $M$
is a Noetherian module, so $M_{s}=M$ for some $s$.  
\end{proof}

\section{Weak equivalences of Hopf algebroids}\label{sec-weak}

In this section, we show that our equivalences of comodule categories
are induced by weak equivalences of Hopf algebroids.  

\begin{definition}\label{defn-weak}
Suppose $\Phi \mathcolon (A,\Gamma)\xrightarrow{}(B,\Sigma)$ is a map of
Hopf algebroids.  We say that $\Phi$ is a \textbf{weak equivalence} if
the induced functor $\Phi_{*}\mathcolon \Gamma \comod
\xrightarrow{}\Sigma \comod$, where $\Phi_{*}M=B\otimes_{A}M$, is an
equivalence of categories.
\end{definition}

We have the following characterization of weak equivalences of Hopf
algebroids. 

\begin{theorem}\label{thm-weak}
A map $\Phi =(\Phi_{0}, \Phi_{1})\mathcolon (A, \Gamma)\xrightarrow{}(B,
\Sigma)$ of flat Hopf algebroids is a weak equivalence if and only if the
composite
\[
A \xrightarrow{\eta_{R}} \Gamma \xrightarrow{\Phi_{0}\otimes 1} B\otimes_{A}\Gamma 
\]
is a faithfully flat ring extension and the map 
\[
B\otimes_{A}\Gamma \otimes_{A}B \xrightarrow{} \Sigma 
\]
that takes $b\otimes x\otimes b'$ to
$\eta_{L}(b)\Phi_{1}(x)\eta_{R}(b')$ is a ring isomorphism.  
\end{theorem}

One can rephrase this theorem using sheaves of groupoids.  A
Hopf algebroid $(A, \Gamma)$ has an associated sheaf of groupoids $\Spec
(A, \Gamma)$ with respect to the flat topology on $\Aff$, the opposite
category of commutative rings (see~\cite{hovey-hopf}).
Hollander~\cite{hollander} has constructed a Quillen model structure on
(pre)sheaves of groupoids in a Grothendieck topology, and
Theorem~\ref{thm-weak} says that $\Phi$ is a weak equivalence if and
only if $\Spec \Phi$ is a weak equivalence of sheaves of groupoids.  

\begin{proof}
The ``if'' half of this theorem is the main result
of~\cite{hovey-hopf}.  Conversely, suppose $\Phi$ is a weak
equivalence.  Then $\Phi_{*}$ is in particular exact, so that $B$ is
Landweber exact over $A$.  Lemma~\ref{lem-land-exact} then guarantees
that $B\otimes_{A}\Gamma$ is flat over $A$.  On the other hand, if
$B\otimes_{A}\Gamma \otimes_{A}M=0$, then $\Phi_{*}(\Gamma
\otimes_{A}M)=0$, so, since $\Phi_{*}$ is an equivalence of
categories, $\Gamma \otimes_{A}M=0$.  Since $A$ is an $A$-module
retract of $\Gamma$, we see that $M=0$.  Hence $B\otimes_{A}\Gamma$
is faithfully flat over $A$.

Now, if $\Phi_{*}$ is an equivalence of categories, then the counit
$\Phi_{*}\Phi^{*}N\xrightarrow{}N$ must be an isomorphism for all
$\Sigma$-comodules $N$.  In particular, $\Phi_{*}\Phi^{*}\Sigma
\xrightarrow{}\Sigma$ must be an isomorphism.  But 
\[
\Phi_{*}\Phi^{*}\Sigma \cong B\otimes_{A} \Phi^{*}(\Sigma \otimes_{B}B)
\cong B\otimes_{A}\Gamma \otimes_{A} B, 
\]
completing the proof. 
\end{proof}

For rings $R$ and $S$, we can have equivalences of categories between
$R$-modules and $S$-modules that are not induced by maps
$R\xrightarrow{}S$; this is, of course, the content of Morita theory.
However, two commutative rings are Morita equivalent if and only if they
are isomorphic.  We view our Hopf algebroids as fundamentally
commutative objects, so we do not expect any non-trivial Morita theory.

\begin{conjecture}\label{conj-Morita}
Suppose $(A, \Gamma)$ and $(B, \Sigma)$ are flat Hopf algebroids such
that the category of $\Gamma$-comodules is equivalent to the category of
$\Sigma$-comodules.  Then $(A, \Gamma)$ and $(B, \Sigma)$ are connected
by a chain of weak equivalences.  
\end{conjecture}

If this conjecture is going to be true, then in particular the
equivalences of comodule categories in Theorem~\ref{thm-main-BP} must be
induced by chains of weak equivalences.  We will now prove this.  

Let $(A, \Gamma)=(BP_{*},BP_{*}BP)$ as usual.  If we have two
Landweber exact $A$-algebras $B$ and $B'$ of the same heights, and a
map $(B,\Gamma_{B})\xrightarrow{}(B',\Gamma _{B'})$ under
$(A,\Gamma)$, then it is a weak equivalence by
Corollary~\ref{cor-maps}.  In general there may be no such map,
though.  We therefore record another obvious source of equivalences.

Suppose we have a groupoid with object set $X$ and morphism set $G$.
Given another set $Y$ and a map $f\mathcolon Y\xrightarrow{}X$, we
previously constructed a groupoid $(Y,G_f)$, where the morphisms in
$G_f$ from $y_1$ to $y_0$ are the morphisms in $G$ from $f(y_1)$ to
$f(y_0)$.  Now suppose we have another map 
$g\mathcolon Y\xrightarrow{}X$, and thus another groupoid $(Y,G_g)$;
we want to know when this is equivalent to $(Y,G_f)$.  Suppose we have
a map $h\mathcolon Y\xrightarrow{}G$ such that 
$\text{target}\circ h=f$ and $\text{source}\circ h=g$, so that $h(y)$
is a morphism from $g(y)$ to $f(y)$ in $G$.  We can then define a
functor $H\mathcolon G_g\xrightarrow{}G_f$ by $H(y)=y$ on objects, and
\[ H(y_0\xleftarrow{u}y_1) = (f(y_0) \xleftarrow{h(y_0)}
           g(y_0) \xleftarrow{u}
           g(y_1) \xleftarrow{h(y_1)^{-1}} f(y_1))
\]
on morphisms (for $u\in G_g(y_1,y_0)=G(g(y_1),g(y_0))$).
Equivalently, let $H'$ be the map
\[ Y\times_{X,g} G\times_{X,g} Y
    \xrightarrow{h\times 1\times(\text{invert}\circ h)} 
    G\times_X G\times_X G \xrightarrow{\text{compose}} G.
\]
Then $H(y_0,u,y_1)=(y_0,H'(y_0,u,y_1),y_1)$.  It is easy to see that
the functor $H$ is an isomorphism of groupoids.

The analogue for Hopf algebroids is as follows.
\begin{lemma}\label{lem-iso-weak-equivalence}
Let $(A,\Gamma)$ be a Hopf algebroid, and suppose $h\mathcolon \Gamma
\xrightarrow{}B$ is a ring homomorphism.  Let $f=h\eta _{L}$ and
$g=h\eta _{R}$.  Then there is an isomorphism of Hopf algebroids from
$(B, \Gamma _{g})$ to $(B,\Gamma_{f})$.  
\end{lemma}

\begin{proof}
 The pair $(A,\Gamma)$ represents a functor from graded rings to
 groupoids, and the conclusion follows from the above discussion by
 Yoneda's lemma.  Alternatively, we can give a formula for the map
 $\Gamma_f\xrightarrow{}\Gamma_g$ as follows.  A map
 $\Gamma_{f}\xrightarrow{}\Gamma _{g}$ of $B$-bimodules is equivalent
 to a map $\Gamma\xrightarrow{}B\otimes_{g}\Gamma\otimes_{g}B$ of
 $A$-bimodules, where the target has the $A$-bimodule structure coming
 from $f$.  This map is the composite
 \[
  \Gamma \xrightarrow{\Delta}
  \Gamma \otimes _{A} \Gamma \xrightarrow{\Delta \otimes 1}
  \Gamma \otimes _{A} \Gamma \otimes_{A} \Gamma 
    \xrightarrow{h\otimes 1\otimes (h\circ\chi)}
  B\otimes _{g} \Gamma \otimes _{g} B  
 \]
 (corresponding to $H'$ in the previous discussion).
\end{proof}



\begin{theorem}\label{thm-stacks}
Let $(A,\Gamma )=(BP_{*},BP_{*}BP)$, and suppose $B$ and $B'$ are
Landweber exact $A$-algebras such that $\height B=\height B'$.  Then the
Hopf algebroids $(B,\Gamma _{B})$ and $(B', \Gamma _{B'})$ are connected
by a chain of weak equivalences.
\end{theorem}

\begin{proof}
Let $C=B\otimes _{A}\Gamma \otimes _{A}B'$.  Let us denote $C$ together
with the ring homomorphism $f\mathcolon A\xrightarrow{}B\xrightarrow{}C$
by $C_{f}$, and $C$ together with the ring homomorphism $g\mathcolon
A\xrightarrow{}B' \xrightarrow{}C$ by $C_{g}$.  Our desired chain of
weak equivalences is
\[
(B, \Gamma _{B}) \xrightarrow{} (C_{f}, \Gamma _{f}) \cong (C_{g},
\Gamma _{g}) \xleftarrow{} (B', \Gamma _{B'}).
\]
The middle isomorphism comes from the evident
map $h\mathcolon \Gamma \xrightarrow{}C$ such that $h\eta _{L}=f$ and
$h\eta _{R}=g$ and Lemma~\ref{lem-iso-weak-equivalence}.  

We now claim that $C_{f}$, and therefore also $C_{g}$, is Landweber
exact.  Indeed, Lemma~\ref{lem-land-exact} implies that $B\otimes
_{A}\Gamma $ and $B'\otimes _{A}\Gamma $ are flat over $A$.  But then
\[
C\otimes _{A}\Gamma =(B\otimes _{A}\Gamma )\otimes _{A}(B'\otimes
_{A}\Gamma )
\]
is also flat over $A$, and so $C_{f}$ is Landweber exact.  

Thus, it suffices to show that $\height C_{f}=\height C_{g}=\height B$.
Because $I_{n}$ is invariant, we have
\[
C/I_{n}\cong (B/I_{n}) \otimes _{A} \Gamma  \otimes _{A}
(B'/I_{n}), 
\]
and therefore $B'/I_{n}=0$ implies $C/I_{n}=0$.  Conversely, suppose
$C/I_{n}=0$, but $B'/I_{n}\neq 0$.  This means that $\height B=\height
B'\geq n$.  Since
\[
B\otimes _{A}(\Gamma \otimes _{A}B'/I_{n})=0,
\]
we conclude that $\Gamma \otimes _{A}B'/I_{n}$ is $v_{\height
B}$-torsion, and therefore $v_{n}$-torsion.  But $B'/I_{n}$ is a
retract of $\Gamma \otimes _{A}B'/I_{n}$ as an $A$-module, so
$B'/I_{n}$ is $v_{n}$-torsion.  Since $B'$ is Landweber exact, this
means $B'/I_{n}=0$, which is a contradiction.  
\end{proof}

\section{The global case}\label{sec-global}

The object of this section is to show that our results about Landweber
exact algebras over $BP_{*}$ extend to Landweber exact algebras over
the complex cobordism ring $MU_{*}$.  Recall that $MU_{*}\cong \Z [x_{1},
x_{2},\dots ]$ for some generators $x_{i}$ of degree $2i$.  All we
require of these generators is that the Chern numbers of $x_{p^{n}-1}$
are all divisible by $p$, as in~\cite{land-exact}.  In this case, the
ideals $I_{p,n}=(p,x_{p-1},\dots ,x_{p^{n-1}-1})$ are invariant and
independent of the choice of generators.  These ideals and $I_{p,\infty
}=\bigcup _{n}I_{p,n}$ are the only invariant prime ideals in
$MU_{*}$~\cite{land-exact}.  

Our first goal is to understand the relation between graded hereditary
torsion theories of $MU_{*}MU$-comodules and graded hereditary torsion
theories of $BP_{*}BP$-comodules.  We use the notation $A_{(p)}$ for
$A\otimes _{\Z}\Zp $, and we recall the well-known fact that
$(MU_{*})_{(p)}$ is a Landweber exact $BP_{*}$-algebra of infinite
height.  Theorem~\ref{thm-main-BP} then gives us an equivalence of
categories between graded $(MU_{*}MU)_{(p)}$-comodules and graded
$BP_{*}BP$-comodules.

\begin{lemma}\label{lem-MU-torsion}
Let $\cat{T}$ be a graded proper hereditary torsion theory of graded
$MU_{*}MU$-comodules, and, for a prime $p$, let $\cat{T}^{(p)}$ denote
the class of $p$-torsion comodules in $\cat{T}$.  Then
$\cat{T}=\bigoplus _{p} \cat{T}^{(p)}$\usc\  that is, $M\in \cat{T}$ if
and only if $M=\bigoplus M_{(p)}$ for $M_{(p)}\in \cat{T}^{(p)}$.
Furthermore, there is a one-to-one correspondence between graded
hereditary torsion theories of graded $p$-torsion $MU_{*}MU$-comodules
and graded proper hereditary torsion theories of graded
$BP_{*}BP$-comodules.
\end{lemma}

Here we refer to a hereditary torsion theory as \textbf{proper} if it
is not the entire category.

\begin{proof}
First of all, if $\cat{T}$ is proper, then $\cat{T}$ must consist
entirely of comodules that are torsion as abelian groups.  Indeed,
suppose $M\in \cat{T}$ is non-torsion.  Let $T(M)$ denote the torsion in
$M$, which is easily seen to be a comodule.  Let $x$ be a nonzero
primitive in $M/T(M)\in \cat{T}$.  Then $x$ is non-torsion.  The
annihilator ideal $I$ of $x$ is invariant, and we claim it is $0$.
Indeed, if $I$ is nonzero, it must
contain a nonzero invariant element of $MU_{*}$, which must be an
integer $m$. But then $mx=0$, contradicting the fact that $x$ is
non-torsion.  The subcomodule of $M/T(M)$ generated by $x$ is thus
isomorphic to $s^{t}MU_{*}$ for some $t$, so $MU_{*}\in \cat{T}$.  This
implies that $\cat{T}$ is the entire category of $MU_{*}$-comodules.
Indeed, we then get $MU_{*}/I\in \cat{T}$ for all invariant ideals $I$.
The Landweber filtration theorem implies that every finitely presented
$MU_{*}MU$-comodule is in $\cat{T}$, and every comodule is a filtered
colimit of finitely presented comodules.  

Now it is easy to check that every torsion comodule $M$ can be written
as $\bigoplus_{(p)} M_{(p)}$, where $M_{(p)}$ is the $p$-localization of
$M$ and therefore is just the $p$-torsion in $M$.  The correspondence
between graded hereditary torsion theories of $p$-torsion
$MU_{*}MU$-comodules and proper graded hereditary torsion theories of
$BP_{*}BP$-comodules follows from the equivalence of categories between
$BP_{*}BP$-comodules and $(MU_{*}MU)_{(p)}$-comodules.  
\end{proof}

We then let $\cat{T}^{(p)}_{n}$ denote the hereditary torsion theory of
$p$-torsion $MU_{*}MU$-comodules corresponding to $\cat{T}_{n}$.  Thus
$\cat{T}^{(p)}_{n}$ is generated by $MU_{*}/I_{p,n+1}$.  For notational
reasons, we let $\cat{T}^{(p)}_{\infty }=(0)$.  

\begin{definition}\label{defn-height-MU}
Given an $MU_{*}$-module $B$ and a prime $p$, define the \textbf{height}
of $B$ at $p$, written $\height_{p}B$, to be the largest $n$ such that
$B/I_{p,n}$ is nonzero, or $\infty$ if $B/I_{p,n}$ is nonzero for all
$n$.  
\end{definition}

We then have the integral analogue of Theorem~\ref{thm-main-BP}. 

\begin{theorem}\label{thm-main-MU}
Let $(A, \Gamma )=(MU_{*},MU_{*}MU)$, and suppose $B$ and $B'$ are two
graded Landweber exact $A$-algebras with $\height _{p}B=\height _{p}B'$ for all
primes $p$.  Then the category of graded $\Gamma _{B}$-comodules is
equivalent to the category of graded $\Gamma _{B'}$-comodules, and both
categories are equivalent to the localization of the category of graded
$\Gamma $-comodules with respect to the torsion theory
$\bigoplus_{p}\cat{T}^{(p)}_{\height _{p}B}$.  This localization is the full
subcategory of graded $\Gamma $-comodules consisting of all those $M$
such that 
\[
\Hom _{A}^{*}(A/I_{p,\height _{p}B+1},M)= \Ext^{1,*}_{\Gamma } (A/I_{p,\height
_{p}B+1},M)=0
\]
for all $p$ such that $\height _{p}B<\infty $.  
\end{theorem}

Note that this theorem implies Theorem~\ref{thm-main-BP}, since if $B$
is a Landweber exact $BP_{*}$-algebra, it is also a Landweber exact
$MU_{*}$-algebra.  

\begin{proof}[Proof of Theorem~\ref{thm-main-MU}]
Theorem~\ref{thm-Giraud} implies that graded $\Gamma
_{B}$-comodules are equivalent to the localization of the category of
graded $\Gamma $-comodules with respect to the kernel $\cat{T}$ of the
functor $M\mapsto B\otimes _{A}M$.  Given a prime $p$, let $\cat{T}^{(p)}$
denote the collection of $p$-torsion comodules in $\cat{T}$.  If $B$ is
zero, there is nothing to prove, so we can assume $B$ is nonzero and
therefore $\cat{T}$ is proper.  Lemma~\ref{lem-MU-torsion} then implies
that we need only check that $\cat{T}^{(p)}=\cat{T}^{(p)}_{\height _{p}B}$.

Suppose first that $\height _{p}B=\infty $.  Then we claim that
$\cat{T}^{(p)}=(0)$.  Indeed, suppose $M$ is a nonzero comodule in
$\cat{T}^{(p)}$.  By choosing a primitive in $M$, we find that
$A/I\in \cat{T}^{(p)}$ for some proper invariant ideal $I$ in
$A$ such that $p^{r}\in I$ for some $r$.  But then $I$ is an
invariant ideal in $A_{(p)}$.  The equivalence of categories
between graded $\Gamma _{(p)}$-comodules and graded
$BP_{*}BP$-comodules preserves invariant ideals.  Since every proper
invariant ideal in $BP_{*}$ is contained in $I_{\infty }$, we see that
$I\subseteq I_{p,\infty }$.  Thus $A/I_{p,\infty }\in \cat{T}$.
Hence $B/I_{p,\infty }=0$.  This means that $1\in I_{p,\infty }B$, so
$1\in I_{p,n}B$ for some $n$.  But then $B/I_{p,n}=0$, violating our
assumption that $\height _{p}B=\infty $.

Now suppose that $\height _{p}B = n <\infty $.  Then $A/I_{p,n}$ is not
in $\cat{T}^{(p)}$ but $A/I_{p,n+1}\in \cat{T}^{(p)}$.
Propositions~\ref{prop-nil} and~\ref{prop-non-nil} imply that
$\cat{T}^{(p)}$ corresponds to the hereditary torsion theory $\cat{T}_{n}$
of $BP_{*}BP$-comodules.   

The characterization of local objects follows from the fact that
$\cat{T}=\bigoplus _{p}\cat{T}^{(p)}$, Lemma~\ref{lem-E-local}, and
Corollary~\ref{cor-main-BP}. 
\end{proof}

We then get analogues of the results of
Sections~\ref{sec-results}--\ref{sec-weak} for $MU_{*}MU$-comodules.  We
will state only the structure theorem for comodules.  We have the same
difficulty with the classification of invariant prime ideals that we
have with invariant radical ideals in the $BP_{*}$-case.  We fix it
analogously.  That is, if $\cat{A}$ is a symmetric monoidal abelian
category with unit $k$ for the tensor product, we define a categorical
ideal $I$ in $\cat{A}$ to be \textbf{categorically prime} if $JK\leq I$
for categorical ideals $J$ and $K$ implies that $J\leq I$ or $K\leq I$.
One checks that Landweber's classification of invariant prime ideals in
$MU_{*}$~\cite{land-ideals} actually classifies categorically prime
invariant ideals.  

\begin{theorem}\label{thm-global-structure}
Let $(A, \Gamma)=(MU_{*},MU_{*}MU)$, and suppose $B$ is a Landweber
exact $A$-algebra.  
\begin{enumerate}
\item [(a)] Every nonzero graded $\Gamma_{B}$-comodule has a nonzero
primitive.  
\item [(b)] The categorically prime invariant ideals in $B$ are
$\{I_{p,n}B \suchthat 0\leq n\leq \height_{p}B \}$.  In particular, this
set contains all the invariant prime ideals.  
\item [(c)] If $B$ is Noetherian, then every graded
$\Gamma_{B}$-comodule that is finitely generated over $B$ admits a
finite filtration by subcomodules 
\[
0 = M_{0} \subseteq M_{1} \subseteq \dots \subseteq M_{s} =M
\]
for some $s$, with $M_{r}/M_{r-1}\cong s^{t}B/I_{p,j}$ for some $s$,$p$,
and $j$ depending on $r$ with
$j\leq \height _{p}B$. 
\end{enumerate}
\end{theorem}

In particular, this theorem applies to $K_{*}K$-comodules or
$\ell_{*}\ell$-comodules, where $K$ is complex $K$-theory, and $\ell$ is
one of the many versions of (periodic, complex oriented) elliptic
cohomology.  

We leave the proof of this theorem to the interested reader, except for
a few comments that illustrate the differences between this theorem and
Theorem~\ref{main-D}.  First of all, in part~(b) we need to assume that
$I$ is categorically prime, rather than just categorically radical.
This is already true in $MU_{*}$, since the ideal $(6)$, for example, is
an invariant radical ideal in $MU_{*}$ not of the desired form.  

Secondly, in the proof of part~(c), we need a Noetherian hypothesis that
is not present in the corresponding fact for $(A,
\Gamma)=(BP_{*},BP_{*}BP)$.  The reason for this is that, in the
$BP_{*}BP$ case, the category of $\Gamma_{B}$-comodules is either
equivalent to the category of $BP_{*}BP$-comodules or to the category of
$E(n)_{*}E(n)$-comodules.  In the first case, we already know the
Landweber filtration theorem, and in the second case $E(n)_{*}$ is
Noetherian.  We believe that Theorem~\ref{thm-global-structure}(c) is
true without the Noetherian hypothesis, however.

\section{$BPJ_{*}BPJ$-comodules}\label{sec-BPJ}

Throughout this section, we let $J$ be a fixed invariant sequence
$p^{i_{0}}, v_{1}^{i_{1}},\dots ,v_{k-1}^{i_{k-1}}$ in $BP_{*}$ of
length $k$.  The spectrum $BPJ$ is constructed from $BP$ by killing this
regular sequence, as in~\cite{johnson-yosimura}, or, in a more modern
fashion, in~\cite[Chapter~V]{elmendorf-kriz-mandell-may}
or~\cite{strickland-products}.  Then $BPJ$ is an associative ring
spectrum, with $BPJ_{*}=BP_{*}/J$.  We will assume that the product on
$BPJ$ has been chosen to be commutative.  This is always possible if
$p>2$, and we believe that it is possible for a cofinal set of ideals
$J$ when $p=2$ although we have not checked the details.  However, it is
not possible when $p=2$ and $J=I_k$.

The co-operation ring $BPJ_{*}BPJ$ is not evenly graded if $k>0$, but
is still free over $BPJ_{*}$, so that $(BPJ_{*},BPJ_{*}BPJ)$ is a Hopf
algebroid.  (When $BPJ$ is not commutative, the structure is more
complicated.)  The object of this section is to extend our results to
Landweber exact $BPJ_{*}$-algebras $B$.  The most important case is
when $J=I_{n}$; the spectrum $BPI_{n}$ is often called $P(n)$.  The
Morava $K$-theory coefficient ring $K(n)_{*}$ is Landweber exact over
$P(n)_{*}$~\cite{yosimura}.

Our first job is to classify the herditary torsion theories of
$BPJ_{*}BPJ$-comodules.  As before, we let $\cat{T}_{n}$
denote the class of all graded $BPJ_{*}BPJ$-comodules that are
$v_{n}$-torsion.  By Lemma~2.3 of~\cite{johnson-yosimura}, $M$ is
$v_{n}$-torsion if and only if $M$ is $I_{n+1}$-torsion.  Of course, any
$BPJ_{*}$-module is automatically $I_{k}$-torsion, so this is only
interesting for $n\geq k-1$.  

The following theorem is our generalization of Theorem~\ref{thm-torsion}.

\begin{theorem}\label{thm-torsion-J}
Let $\cat{T}$ be a graded hereditary torsion theory of graded
$BPJ_{*}BPJ$-comodules, and suppose that $\cat{T}$ contains some nonzero
comodule that is finitely presented over $BPJ_{*}$.  Then
$\cat{T}=\cat{T}_{n}$ for some $n\geq k-1$.  
\end{theorem}

This theorem is proved just as Theorem~\ref{thm-torsion}, except that
the results of~\cite{landweber-comodules}, which we used in the proof of
Lemma~\ref{lem-nil}, are not written so as to apply to
$BPJ_{*}BPJ$-comodules.  So one can either reprove the results
of~\cite{landweber-comodules} in this case, or construct a direct proof
of Lemma~\ref{lem-nil} for $BPJ_{*}BPJ$-comodules using the results
of~\cite{johnson-yosimura}.  

We can then define the height for $BPJ_{*}$-algebras.

\begin{definition}\label{defn-height-J}
Suppose $B$ is a nonzero graded $BPJ_{*}$-module.  We define the
\textbf{height} of $B$, written $\height B$, to be the largest $n$ such
that $B/I_{n}$ is nonzero, or $\infty$ if $B/I_{n}$ is nonzero for all
$n$.  
\end{definition}

Note that every nonzero $BPJ_{*}$-module $B$ has $\height B\geq k$.

Here is the analogue of Theorem~\ref{thm-main-BP}.  

\begin{theorem}\label{thm-main-BPJ}
Let $(A, \Gamma )=(BPJ_{*}, BPJ_{*}BPJ)$, and suppose $B$ and $B'$ are
two graded Landweber exact $BPJ_{*}$-algebras with $k\leq \height
B=\height B'=n \leq \infty $.  Then the category of graded $\Gamma
_{B}$-comodules is equivalent to the category of graded $\Gamma
_{B'}$-comodules.  If $n=\infty$, these categories are equivalent to the
category of graded $\Gamma$-comodules.  If $n<\infty$, these categories
are equivalent to the localization of the category of graded
$\Gamma$-comodules with respect to the torsion theory $\cat{T}_{n}$.
\end{theorem}

We then get analogues of the results of
Sections~\ref{sec-results}--\ref{sec-weak} for $BPJ_{*}BPJ$-comodules.
These depend on the results of~\cite{johnson-yosimura} on the structure
of $BPJ_{*}BPJ$-comodules to replace the results of Landweber on the
structure of $BP_{*}BP$-comodules.  

In particular, we get versions of the Miller-Ravenel, Morava, and
Hovey-Sadofsky change of rings theorems.  These use the spectra
$v_{n}^{-1}BPJ$ and $E(n,J)$ for $n\geq k$ in place of
$v_{n}^{-1}BP$ and $E(n)$.  Here $E(n,J)_{*}$ is
Landweber exact over $BPJ_{*}$ with 
\[
E(n,J)_{*}\cong v_{n}^{-1}(BPJ_{*}/(v_{n+1}, v_{n+2},\dotsc )).
\]

Here is the structure theorem for comodules.  

\begin{theorem}\label{thm-structure-J}
Let $(A, \Gamma)=(BPJ_{*},BPJ_{*}BPJ)$, and suppose $B$ is a Landweber
exact $A$-algebra.  
\begin{enumerate}
\item [(a)] Every nonzero graded $\Gamma_{B}$-comodule has a nonzero
primitive.  
\item [(b)] The categorically radical invariant ideals in $B$ are
$\{I_{n}B \suchthat k\leq n \leq \height B\}$. 
\item [(c)] Every graded $\Gamma_{B}$-comodule that is finitely
presented over $B$ admits a finite filtration by subcomodules
\[
0 = M_{0} \subseteq M_{1} \subseteq \dots \subseteq M_{s} =M
\]
for some $s$, with $M_{r}/M_{r-1}\cong s^{t}B/I_{j}$ for some $s$,$p$,
and $j$ depending on $r$ with $k\leq j\leq \height B$.
\end{enumerate}
\end{theorem}

\providecommand{\bysame}{\leavevmode\hbox to3em{\hrulefill}\thinspace}
\providecommand{\MR}{\relax\ifhmode\unskip\space\fi MR }
\providecommand{\MRhref}[2]{%
  \href{http://www.ams.org/mathscinet-getitem?mr=#1}{#2}
}
\providecommand{\href}[2]{#2}

\end{document}